\documentclass[11pt, a4paper, reqno, oneside]{amsart} 
\usepackage[foot]{amsaddr}
\usepackage[utf8]{inputenc}
\usepackage[UKenglish]{babel} 
\usepackage[a4paper]{geometry}
\usepackage{graphicx}
\usepackage{quiver}
\usepackage{amsmath,amssymb,amsfonts,amsthm}
\usepackage{hyperref}

\usepackage{heuristica}
\usepackage[heuristica,vvarbb,bigdelims]{newtxmath}
\usepackage[T1]{fontenc}

\usepackage{comment}
\usepackage[final]{showlabels} 
\usepackage[]{todonotes} 

\usepackage[style = alphabetic, backend = bibtex, abbreviate = true, arxiv = abs, doi = false, isbn = false, url = false, giveninits = true]{biblatex} 
\renewbibmacro{in:}{}

\theoremstyle{plain}
\newtheorem{theorem}{Theorem}
\newtheorem{proposition}{Proposition}[section]
\newtheorem{lemma}[proposition]{Lemma}
\newtheorem{corollary}[proposition]{Corollary}

\theoremstyle{definition}
\newtheorem{definition}{Definition}
\newtheorem*{definition*}{Definition}

\newtheorem*{hypothesis}{Hypothesis}

\theoremstyle{remark}
\newtheorem{example}{Example}
\newtheorem{remark}[proposition]{Remark}

\numberwithin{equation}{section}

\newcommand{\Z}{\ensuremath{\mathbb{Z}}}
\newcommand{\R}{\ensuremath{\mathbb{R}}}
\newcommand{\CC}{\ensuremath{\mathbb{C}}}

\newcommand{\Q}{\ensuremath{\mathbb{Q}}}
\newcommand{\N}{\ensuremath{\mathbb{N}}}

\newcommand{\uhp}{\ensuremath{\mathbb{H}}}

\newcommand{\norm}[1]{\left\lVert#1\right\rVert}

\DeclareMathOperator{\diag}{diag}
\DeclareMathOperator{\vol}{vol}
\DeclareMathOperator{\SO}{SO}
\DeclareMathOperator{\SL}{SL}
\DeclareMathOperator{\GL}{GL}
\DeclareMathOperator{\PGL}{PGL}
\DeclareMathOperator{\stab}{Stab}
\DeclareMathOperator{\id}{id}

\usepackage{microtype}

\title[Sup-norms in higher rank, level aspect]{The sup-norm problem for newforms of large level on $\PGL(n)$}
\author[Radu Toma]{Radu Toma}
\address{Sorbonne Université, Université Paris Cité, CNRS, IMJ-PRG, F-75005 Paris, France}
\email{toma@imj-prg.fr}
\thanks{The author was supported by the Germany Excellence Strategy grant EXC-2047/1-390685813 through the Bonn International Graduate School of Mathematics and by the European Research Council Advanced Grant 101054336.}

\begin{document}
	
\begin{abstract}
	Let $N$ be a prime and $\phi$ be a Hecke-Maaß cuspidal newform for the Hecke congruence subgroup $\Gamma_0(N)$ in $\SL_n(\R)$. We define a notion of the bulk $\Omega_N$ of the space $\Gamma_0(N) \backslash \SL_n(\R)/ \SO(n)$ with respect to a compact set $\Omega \subset \SL_n(\R)$. For any prime $n$, we prove sub-baseline bounds for the sup-norm of $\phi$ restricted to $\Omega_N$. Conditionally on GRH, we generalise this result to all $n \geq 2$. The methods involve a new reduction theory with level structure, based on generalisations of Atkin-Lehner operators.
\end{abstract}

\maketitle

\section{Introduction}

Let $n \geq 2$ be an integer. 
This article is concerned with bounding the sup-norm of Hecke-Maaß forms on the space
\begin{displaymath}
X_n(N) = \Gamma_0(N) \backslash \SL_n(\R) / \SO(n)
\end{displaymath} 
in terms of the parameter $N$, called the \emph{level}. 
Here, $\Gamma_0(N) \leq \SL_n(\Z)$ is the subgroup of integral matrices with last row congruent to $(0, \ldots, 0, \ast)$ modulo $N$, where $\ast$ stands for any non-zero residue class.

We normalise the invariant measure on $X_n(N)$ so that it has volume asymptotically equivalent to $N^{(n-1)}$. 
Now let $\phi$ be a Hecke-Maaß form on this space, that is, a square-integrable joint eigenfunction of the invariant differential operators and the unramified Hecke algebra. 
Assuming that $\norm{\phi}_2 = 1$, the sup-norm problem asks for non-trivial bounds on $\norm{\phi}_\infty$. 
Several parameters can be considered for this question, the most studied being the spectral parameter and the level.

\subsection{Some history}
This problem has a rich history and the first breakthrough in the eigenvalue aspect for $n=2$ was achieved by Iwaniec and Sarnak \cite{IS}. 
They prove that $\norm{\phi}_\infty \ll_{N, \varepsilon} \lambda^{5/24 + \varepsilon}$ for any $\varepsilon > 0$. 
This is an improvement over the so-called local bound $\norm{\phi}_\infty \ll_{N} \lambda^{1/4}$. 
Their method of using an amplified pre-trace formula remains one of the main tools for obtaining such non-trivial, sub-local bounds.

In the level aspect, the baseline bound expected to hold is $\norm{\phi}_\infty \ll_{\lambda, \varepsilon} N^{\varepsilon}$ for $\phi$ a newform.
The first improvement for $n=2$ is due to Blomer and Holowinsky \cite{BHolo}, with important refinements by Harcos and Templier \cite{HT2}, \cite{HT3}, and the current record bound $\norm{\phi}_\infty \ll_{\lambda, \varepsilon} N^{1/4 + \varepsilon}$ is due to Khayutin, Nelson and Steiner \cite{KNS}.
These papers deal with the case of square-free level $N$, and bounds for general $N$ were achieved in \cite{saha}. 
The fact that much of the work on this problem historically focused on square-free levels is in large part a consequence of using Atkin-Lehner operators. 
This aspect of the problem forms one of the main topics of this paper.

Though many other variations of the problem exist, we consider now its development in higher rank, that is, for $n > 2$. 
In the spectral aspect we only mention here a selection, namely the work of Blomer and Pohl \cite{BP} (for $\operatorname{Sp}_4$), Blomer and Maga \cite{BMn} (for $\SL_n$), and Marshall \cite{marshall} (for more general Lie groups). 
They achieve power savings over the local bound for any $n \geq 2$, though they only consider the sup-norm of automorphic forms restricted to a fixed compact set. 
The implicit constants in their bounds thus depend on this set, which is necessary by a result of Brumley--Templier \cite{BT} about large values in the cusp. 
An investigation of the global sup-norm is the topic of Blomer, Harcos and Maga's paper \cite{BHMn}.

The present article deals with the sup-norm problem in higher rank, in the level aspect. 
Despite the progress described above, there are very few results in this setting. 
The first result, due to Hu \cite{hu}, considers the case of prime-power levels $N = p^c$, where $c$ is large, with $\phi$ corresponding to a so-called minimal vector, thus not applying to newforms. 
These forms are more suitable for the $p$-adic methods employed by Hu. 
Similar to many results in the spectral aspect, the bounds are given for the sup-norm of the restriction to a fixed adelic compact set.

The second result \cite{toma} is due to the author of this paper and concerns automorphic forms on a different family of spaces $\Gamma \backslash \SL_n(\R) / \SO(n)$, where $\Gamma$ is a subgroup coming from the units of an order in a division algebra of degree $n$. 
These spaces are compact, and the bounds provided are global and in terms of their volume. 
The degree $n$ is restricted to prime numbers and results can only be extended partially to odd degrees.

Moreover, the argument is based on the fact that proper subalgebras of division algebras of prime degree are automatically fields, and that zero is the only element of norm zero.
The situation is decidedly different for the matrix algebra, whose orders give rise to the groups $\Gamma_0(N)$, and thus the methods of \cite{toma} seem to be insufficient in this case.

Not only throughout the history of the sup-norm problem, but also of the subconvexity problem, the level aspect, particularly for prime or square-free levels, is often the last one to be successfully tackled.
Given its significance in number theory, this suggests a serious, general difficulty and a need for new ideas.

\subsection{Statement of results} \label{sec-results}

In this paper, we consider Hecke-Maaß cuspidal newforms on $X_n(N)$ for $n \geq 2$ and $N$ prime. 
Motivated by \cite{BT} and the necessity of doing so in the spectral aspect, we truncate cuspidal regions and focus on the ``bulk'' of $X_n(N)$, which we construct as follows.

Let $\Omega \subset X_n(1)$ be a fixed compact set and define the preliminary set $\Omega_N^\natural \subset X_n(N)$ as the preimage of $\Omega$ under the natural projection $X_n(N) \to X_n(1)$. 
It is easy to see that $\vol(\Omega_N^\natural) \asymp_{\Omega} \vol(X_n(N)) \sim N^{n-1}$.
However, some components of the preimage spread out into the cusps of $X_n(N)$ and need to be truncated.
We use a reduction theory with level to describe this, inspired by the case $n=2$.

For this, we study a new generalisation of the classical Fricke involution on $X_n(N)$, denoted by $W_N$.
We show that $\Omega_N^\natural$ can be covered by $S \cup W_N S$ for a Siegel set $S$ given in terms of Iwasawa coordinates by
\begin{displaymath}
	\{ z = n(x)a(y) : |x_{ij}| \leq 1/2 \text{ for }1 \leq i < j \leq n \text{ and } y_1 \gg 1/N, y_i \gg 1 \text{ for } 2 \leq i \leq n-1 \}.
\end{displaymath}
Most of the volume of $S$ is at $y_1 \asymp 1/N$ and $y_i \asymp 1$ for all $i \geq 2$.

We now further truncate $\Omega_N^\natural$ by removing a cuspidal region $C_N$ with significantly smaller volume
\begin{displaymath}
	\vol(C_N) \ll N^{(n-1)(1-1/2n^4)}
\end{displaymath} 
that is given in the domain above by the additional bounds 
\begin{equation} \label{eq:exceptional}
	y_1 \gg N^{-1 + 1/2n^4} \text{ or } y_i \gg N^{1/2n^4} \text{ for some } i \geq 2.
\end{equation}
Define now the \emph{bulk} as 
\begin{equation} \label{eq:bulk} 
	\Omega_N = \Omega_N^\natural \setminus C_N,
\end{equation}
which, volumewise, constitutes an arbitrarily large  proportion of the space $X_n(N)$ for large enough $N$ and $\Omega$.
For more details, see Section \ref{sec-red}.

We prove two new results, the first of which applying to all \emph{prime} $n \geq 2$. 

\begin{theorem} \label{thm-uncond}
	Let $n$ and $N$ be primes. Let $\phi$ be a Hecke-Maaß cuspidal newform on $X_n(N)$ with spectral parameter $\mu$ and define $\Omega_N \subset X_n(N)$ as in \eqref{eq:bulk} with respect to a fixed compact set~$\Omega \subset X_n(1)$. 
	For large $N$, we have the bound
	\begin{displaymath}
		\norm{\phi |_{\Omega_N}}_\infty \ll_{\Omega, n, \mu, \varepsilon} N^{-\frac{1}{4 n^2} + \varepsilon}.
	\end{displaymath}
\end{theorem}

The proof involves understanding the geometric structure of the problem as well as handling rather delicate diophantine conditions. 
It is the latter that are not yet well enough understood in the case where $n$ is not prime.

Additionally, as in the case of $n=2$, bounding $\phi$ away from the bulk, e.g. on $C_N$, requires strong bounds from the Whittaker expansion of $\phi$.
In the level aspect, these are not available in higher rank, and further investigations present serious obstacles.
These remain beyond the scope of the present paper (however, see the Appendix for some remarks). 

However, the geometric ideas introduced here are valid in full generality and already capture a significant part of the problem. 
To support this claim, we present below results for all $n \geq 2$, even improving those above numerically, assuming the existence of an efficient amplifier.
 
For this, let $\lambda(p)$ be the Hecke eigenvalue of $\phi$ for the Hecke operator $T_p$, where $p$ is a prime not dividing $N$, normalised so that $\lambda(p) \ll p^{(n-1)/2}$ under the Ramanujan-Petersson conjecture. See Section \ref{sec-hecke-alg} for a precise definition.
\begin{hypothesis}
	Let $\delta > 0$ be any positive constant and $N \gg_{\delta, \mu} 1$ be large enough. If $L \gg N^\delta$, then
	\begin{equation} \label{eq-hypothesis}
	\sum_{p \in \mathcal{P}} \frac{|\lambda(p)|}{p^{(n-1)/2}} \gg_{\varepsilon} L^{3/4-\varepsilon}.
	\end{equation}
\end{hypothesis}
We prove in Lemma \ref{lemma-hypothesis} that condition \eqref{eq-hypothesis} is true assuming the Grand Riemann Hypothesis. 
It is similar to condition (1.24) in \cite{IS}, which is checked in \cite{huang} for dihedral Maaß forms and in \cite{young} for Eisenstein series and leads to an improved exponent in the bound of Iwaniec and Sarnak, as explained in \cite[Remark 1.6]{IS}.

\begin{theorem} \label{thm-main}
	Let $n \geq 2$ and $N$ be a prime. Let $\phi$ be a Hecke-Maaß cuspidal newform on $X_n(N)$ with spectral parameter $\mu$ and define $\Omega_N \subset X_n(N)$ as in \eqref{eq:bulk} with respect to a fixed compact set~$\Omega$. Assuming hypothesis \eqref{eq-hypothesis}, we have the bound
	\begin{displaymath}
		\norm{\phi |_{\Omega_N}}_\infty \ll_{\Omega, n, \mu, \varepsilon} N^{-\frac{1}{8n} + \varepsilon}.
	\end{displaymath}
	In particular, the bound holds under the Grand Riemann Hypothesis.
\end{theorem}


Considering previous work on the sup-norm problem in higher rank, the main contribution of this paper is a new counting argument, based on the reduction of the domain $\Omega_N$ using generalised Atkin-Lehner operators, which might be of independent interest. 
These arguments significantly generalise and give a new perspective on the geometric methods of Harcos and Templier \cite{HT3}, which generated many strong results for the sup-norm problem on $\GL(2)$ (e.g.\ \cite{BHMM}, \cite{saha}, \cite{assing}).
They also seem to be fundamentally different and provide stronger results than in the spectral aspect in higher rank, where savings are inverse super-exponential in $n$ \cite{gillman}, as opposed to our inverse polynomial savings.
In any case, the methods presented here provide the first steps in tackling the level aspect in higher rank and, we believe, a useful framework for proving more general and possibly stronger results in the future.

\subsection{Methods}

For proving both main theorems, we employ an amplified pre-trace formula to transform the analytic issue of bounding the sup-norm into a counting problem. 
This is one of the most common methods of studying the sup-norm of automorphic forms and goes back to the influential paper \cite{IS}. 

As in Proposition \ref{prop-amp} below, we reduce the problem of bounding $\phi(z)$ for $z \in \SL_n(\R)$ to counting matrices in sets of the form
\begin{displaymath}
	H(z, m, N) := \{ \gamma \in \mathcal{M}_n(\Z,N) \mid \det \gamma = m,\, z^{-1} \gamma z = O( m^{1/n}) \},
\end{displaymath}
where $m$ is running over different, potentially sparse, sets of integers. 
Here, $\mathcal{M}_n(\Z,N)$ is the set of integral matrices with last row congruent to $(0, \ldots, 0,*)$ modulo $N$.
This is an order in the algebra of rational matrices.

\subsubsection{Lattices} \label{sec-intro-lats}

To start, we give the sets $H(z, m, N)$ an interpretation in terms of lattices, which motivates the development of new tools introduced below.
This is natural, since we recall that the space $X_n(1)$ parametrises shapes of unimodular lattices by associating to $z \in \SL_n(\R)$ the lattice
\begin{displaymath}
	L = \Z^n \cdot z \subset \R^n.
\end{displaymath}
Here we understand $\R^n$ and $\Z^n$ as sets of row vectors.
In this interpretation, the matrix $z$ gives a specific basis for $L$.
If $N$ is prime, the space $X_n(N)$ now parametrises \emph{pairs} $(L, L_N)$ of lattices, up to simultaneous rotation by $\SO(n)$, where
\begin{displaymath}
	L_N = \Z^n \cdot \diag(N, \ldots, N, 1) z = (N\Z \times \cdots N\Z \times \Z ) \cdot z,
\end{displaymath}
is a sublattice of $L$.

Let $e_1, \ldots, e_n$ be the standard basis for $\R^n$. We evaluate the condition
\begin{displaymath}
	z^{-1} \gamma z = O(m^{1/n})
\end{displaymath}
at the vectors $e_i$, after multiplying from the left by $z$.
This amounts to the conditions
\begin{displaymath}
	e_i \cdot \gamma z \in B(m^{1/n} \norm{e_i \cdot z})
\end{displaymath}
for each $i$, where $B(r)$ is a Euclidean ball of radius $O(r)$ around $0$.
Note that, since $\gamma$ is an integral matrix, $e_i \cdot \gamma z$ is a lattice point in $L$ determining the $i$-th row of $\gamma$.
Moreover, it is important to observe that $e_n \cdot \gamma z$ is additionally a lattice point in the sublattice $L_N$.
On the other side, $e_i \cdot z$ is simply one of the basis vectors in the basis of $L$ determined by $z$.

To count the number of relevant $\gamma$, we can therefore bound the number of possibilities for each of their rows and by the conditions above we reduce to counting lattice points in balls. 
However, this naïve strategy needs to be refined by an application of the Gram-Schmidt process, which we make precise in Section \ref{sec-counting-strategy}. By its very nature, this involves the Iwasawa coordinates of $z$.

In any case, it is apparent that the dependence on $z$ manifests itself in two ways already at this level.
Firstly, there might be many lattice points that we count because the basis vectors $e_i \cdot z$ which control the size of the balls are large.
Secondly, the lattices $L$ and $L_N$ might be very dense, in the sense that they could have very short vectors relative to their covolume.

Understanding such issues is one of the main goals of reduction theory and the geometry of numbers. 
However, the level structure needs to be taken into consideration and, indeed, puts serious restrictions on the prospect of success for the amplified pre-trace formula strategy.
We develop a novel reduction theory with level structure in Section \ref{sec-red}, and we describe the main ideas below.

\subsubsection{Generalised Atkin-Lehner operators and reduction}

In a nutshell, classical reduction theory provides a way to fit a fundamental domain for $X_n(1)$ inside a Siegel set (for the cusp at infinity). 
If $z \in \SL_n(\R)$ lies in such a fundamental domain, its rows then provide a reduced basis for the lattice $L$, that is, a basis of vectors that are as short and as orthogonal as possible.
We also obtain in this way an interpretation of the Iwasawa coordinates of $z$ in terms of the successive minima of $L$. 
This is not only important for implementing the refined counting strategy described above, but also for compensating with other tools when the latter fails.

For instance, solving the matrix counting problem optimally and plugging the result into the amplified pre-trace formula \emph{cannot} yield sub-baseline bounds when $z$ is high enough in the cusp.
One then compensates by using the Fourier expansion, which gives strong bounds in terms of the Iwasawa $y$-coordinates following from the cuspidality of our automorphic form $\phi$.
This is common for many of the previous works \cite[Lemma A.1]{IS}, \cite[Lemma 5.1]{HT2}, etc.

In the level aspect, already $z = \id_n$ has to be treated using the Fourier bound.
Such a point lies in $\Omega_N^\natural$, for $\Omega$ a compact neighbourhood of the identity, and can be considered as lying in a cuspidal region for this problem.
It is part of the set $C_N$ removed from $\Omega_N^\natural$ to form the bulk $\Omega_N$ in \eqref{eq:bulk}.

From one perspective, which we do not explicate here further, the difficulty here is caused by the contribution of Eisenstein series on the spectral side of the pre-trace formula.
In our framework, the reason is that, even though $L$ is a perfectly balanced lattice and $z$ gives an actual orthogonal basis of shortest vectors, the sublattice $L_N$ is maximally imbalanced.

A desirable reduction theory with level structure might thus fulfil the following.
It should provide a basis for the lattice $L$ that, while perhaps not reduced, gives useful information about the shortest vectors in the sublattice $L_N$ and about the Iwasawa coordinates, meaning the Gram-Schmidt process for the basis.
It should also permit some understanding of the successive minima of both $L$ and $L_N$.
Of course, preserving the level structure means changing bases is only allowed by matrices in $\Gamma_0(N)$.
However, there are additional symmetries at our disposal.

It was recognised early on in the treatment of the sup-norm problem in the level aspect that Atkin-Lehner operators would be useful for such reductions.
It is classically not hard to see that one can fit the fundamental domain for $X_2(N)$ with respect to the action of these operators in a single Siegel set of \emph{finite} volume for $N$ square-free.
This is because the Atkin-Lehner operators for $N$ square-free conjugate all cusps to the cusp at infinity.
Unfortunately, for powerful levels there is a deficiency of Atkin-Lehner operators and this forms an important reason why the first and many results on the sup-norm problem are restricted to square-free levels.

The first authors to connect these group theoretic facts to lattices were Harcos and Templier in \cite[Lemma 2.2]{HT2}.
For example, at the level of lattices, the Fricke involution for prime levels can be understood as switching the lattices in the pair $(L, L_N)$.
Together with ideas from reduction theory, this allows us to trade imbalancedness of $L$ or $L_N$ for closeness of $z$ to the cusp (see loc.\ cit.).
Effectively, when the matrix counting results are weak, the Fourier bound should be better. 

Generalising the case $n = 2$, we study the symmetries of $X_n(N)$. 
The point of departure from the classical case is the observation that $\PGL(n)$ for $n > 2$ has an additional outer automorphism, given by $z \mapsto z^{-T}$. 
This corresponds to taking duals, either at the level of lattices, or at the level of automorphic forms.
In this paper, we use this to introduce in Section \ref{sec-fricke} a higher-rank Atkin-Lehner operator corresponding to the \emph{Fricke involution}.
It has probably been implicitly present in the theory of newforms, yet an explicit definition seems hard to find in the literature.
\begin{definition*}
	Let
	\begin{displaymath}
		A_N = N^{-1/n} \diag(1, \ldots, 1, N)
	\end{displaymath} 
	and define the \emph{Fricke involution} $W_N: L^2(X_n(N)) \longrightarrow L^2(X_n(N))$ as
	\begin{displaymath}
		W_N \phi(z) = \phi(A_N \cdot z^{-T}).
	\end{displaymath}
\end{definition*}

We also perform an investigation of other potential generalisations of Atkin-Lehner operators. 
First, we prove that the normaliser of $\Gamma_0(N)$ in $\operatorname{PGL}_n(\R)$, the source of Atkin-Lehner operators for $n=2$, is \emph{trivial} for $n > 3$. 
We refer to Section \ref{sec-normaliser}.

\begin{theorem} \label{thm-normaliser-R}
	For $n > 2$, the normaliser of $\Gamma_0 (N)$ inside $\PGL_n(\R)$ is trivial.
\end{theorem}

We then provide a different perspective on the classical Atkin-Lehner operators and show in Corollary \ref{prop-only-fricke} that the only possible generalisation in this interpretation is the Fricke involution. 
On the one hand, this is in contrast to the case of square-free levels in $\PGL(2)$, but it is also a reflection of the remarkable \emph{lack} of such symmetries for powerful levels. 
Therefore, we first only consider the case of prime level in this paper, similar to the common restrictions in the rank-one case.

The main result of our reduction theory is given in Proposition \ref{prop-fricke-red}. 
It satisfies the intuition from the $n = 2$ case, where the bulk of the reduced fundamental domain is at $\Im(z) \asymp 1/N$.
In general, there are the Iwasawa coordinates $y_1, \ldots, y_{n-1}$ and the bulk (in terms of volume) can be found at 
\begin{displaymath}
	y_1 \asymp \frac{1}{N},\, y_2 \asymp \ldots \asymp y_{n-1} \asymp 1.
\end{displaymath}
More generally, one can reduce $\Omega_N^\natural$ using the Fricke involution to a subset where $y_1 \gg 1/N$ and $y_i \gg 1$ for $i \geq 2$.
Roughly speaking, if any of these coordinates are significantly larger than in the bulk, the counting results become too weak, as is to be expected (see the previous paragraphs and the work of Harcos--Templier).
Ideally, a Fourier bound would suffice to treat these cases, yet this tool is not currently available in the necessary form.
For this reason, we excise the cuspidal set $C_N$ in the statements of the main theorems.
As noted previously, this set has significantly smaller volume than $\Omega_N^\natural$.

However, the reduction of the full fundamental domain for $\Gamma_0(N)$ is more complex, as can be seen from the case work in Section \ref{sec-fricke-red}.
It seems that more refined information can be extracted and doing so would be an important next step in the study of the sup-norm problem in the level aspect.

In higher rank, the reduction process involves the outer automorphism included in the Fricke involution and thus dualising lattices.
We are therefore required to develop tools for keeping track of sizes of vectors in the lattices associated to $z$ and its conjugate under the Fricke involution, as well as their duals.
This is the content of Section \ref{sec-two-lattices} and Table \ref{table}.
We have found the language of wedge products particularly useful for this because of its flexibility in relating lengths of vectors in lattices and their duals with Iwasawa coordinates.

As a historical interlude, we point out some connections of the above considerations with previous work.
The Atkin-Lehner involutions were already used in the breakthrough \cite{BHolo}, but balancedness of lattices was interpreted in terms of Diophantine approximation properties of the Iwasawa coordinates, using terminology from the circle method.

The language of lattices was used directly in \cite{HT2}, \cite{HT3}, and subsequent works, and lead to strong numerical improvements to the bounds. 
However, the counting problem is interpreted using coordinates not truly inherent to lattices. 
Many computations in the $\GL(2)$ case use, in fact, the ``sporadic'' symplectic nature of this group.
This is not available in higher degree and the direct use of coordinates seems to be very cumbersome.

For the family of groups $\PGL(n)$, some ideas reminiscent of the more general strategy used here can be seen in \cite[Sec.\ 3.2]{BHMn}.
We refer also to \cite{venk}, where certain aspects of the geometry of $X_n(N)$ are studied using lattices as well.

\subsubsection{Detecting sparse sequences of determinants}

The upshot of the reduction theory with level and the iterative counting strategy is that we get bounds for the set
\begin{displaymath}
	\bigcup_{1 \leq m \leq \Lambda} H(z, m, N)
\end{displaymath} 
for a parameter $\Lambda$ small enough in terms of $N$, uniformly in the balanced part of $\Omega_N$.
The motto of the counting strategy under these conditions is a rigidity principle: \emph{the last row of $\gamma \in H(z, m, N)$ determines the whole matrix}.

However, the unconditional amplifier of \cite{BM4} gives rise to a counting problem where matrices have perfect power determinants, for instance, $n$-th powers. 
Such a sequence of determinants is too sparse and the method above, averaging over all determinants, produces gross over-counting.
Similar issues are well-known already in the classical case $n = 2$ (see e.g.\ the special treatment of square determinants in \cite{HT3}).

The appearance of sparse sequences of determinants on the geometric side is due to the lack of good lower bounds for Hecke eigenvalues. 
Indeed, such bounds are precisely what Hypothesis \eqref{eq-hypothesis} provides. 
Unconditionally, there is thankfully a substitute obtained from Hecke relations, such as $\lambda(p)^2 - \lambda(p^2) = 1$ in suitable normalisation for $n = 2$, from which one derives that at least one of the two eigenvalues is bounded from below.
Introducing the Hecke operator $T_{p^2}$ in this way results in sequences of square determinants, and we have similar phenomena in higher degree.

We are able to detect perfect power determinants by using a refinement of the counting strategy above (see Section \ref{sec-power-dets}). 
The problem reduces to counting solutions to an equation of the shape
\begin{displaymath}
	\chi_\gamma(X) - Y^\nu = 0
\end{displaymath}
for $1 \leq \nu \leq n$, where $\chi_\gamma$ is the characteristic polynomial of $\gamma$.
If this equation is irreducible, then a powerful theorem of Heath-Brown \cite{HB} provides an adequate non-trivial bound.  

To treat the case where the polynomial is reducible, we assume that $n$ is prime to simplify the classification of these degenerate cases.
We can thus reduce to counting matrices with $\chi_\gamma(X) = (X - m)^n$.
For $n = 2$, this is the special case of parabolic matrices that was also handled in \cite[Lemma 4.1]{HT2}.

Finally, resolving this problem involves some group theoretic investigations once more. 
We classify the cusps of $X_n(N)$ as in Lemma \ref{lemma-cusps}, of which there are $n$ many, and observe the action of the Fricke involution on them.
The cusp corresponding to the identity element, informally the cusp at infinity, can be dealt with by the counting methods already introduced.
The one corresponding to the long Weyl element is conjugated to the identity by the Fricke involution.

Counting at ``intermediate'' cusps presents new challenges, which might be a consequence of the lack of more symmetries of $X_n(N)$ for $n > 2$.
Although much of what is developed in this paper appears to the author to be conceptually necessary and inherent to the problem, this last step is solved by a trick, as one might call it.
We use the specific shape of the amplifier of Blomer and Maga.
Namely, we take advantage of the fact that, for certain Hecke sets attached to primes $p$ and $q$, the determinantal divisors are asymmetric in terms of $p$ and $q$, as in \eqref{eq-hecke-sets}.
This eventually collapses an average over two primes to one over a single prime (the case $p = q$), and leads to the required power saving.

\subsection*{Notation}
By the Vinogradov notation $f(x) \ll g(x)$ for two functions $f, g$ it is meant that $|f(x)| \leq C \cdot |g(x)|$, at least for large enough $x$, for some $C > 0$ called the implied constant.
Similarly, for a matrix $X$ and a scalar function $f(X)$ we say that $X = O(f(X))$ when $\norm{X} \leq C \cdot f(X)$ for some constant $C > 0$ and some choice of matrix norm $\norm{\cdot}$.

We use $\ll_P$ to say that the implied constant depends on a parameter $P$, yet we do not always add the subscript if it is clear from context in order to avoid clutter.
For instance, dependency on the compact space $\Omega \subset \SL_n(\R)$ includes dependency on $n$. 

\section{Preliminaries on lattices} \label{sec-prelims}

Consider the real vector space $V = \R^n$ with standard inner product $\langle v, w \rangle = v \cdot w^T$, where we think of  $v, w \in V$ as \emph{row} vectors in the standard basis $e_1, \ldots, e_n$. Let $z$ be a matrix in $\GL_n(\R)$ and define $L_z$ to be the lattice $\Z^n \cdot z$ inside $V$. Note that $e_i \cdot z$ is equal to the $i$-th row of $z$. We also define the inner product and norm
\begin{displaymath}
    \langle v, w \rangle_z = \langle vz, wz \rangle, \quad \norm{v}_z = \sqrt{\langle vz, vz \rangle},
\end{displaymath}
for $v, w \in V$.

The \emph{dual} lattice $L_z^\ast$ is defined as the set of vectors $w$ such that $\langle v, w \rangle \in \Z$ for all $v \in L_z$. It is straight-forward to compute that
\begin{displaymath}
    L_z^\ast = L_{z^{-T}}.
\end{displaymath}
We also note that $L_z = L_w$ for any $w \in \GL_n(\Z) \cdot z$.

\subsection{Exterior powers}
	
If $k$ is a positive integer, the $k$-th exterior power of $L_z$ is denoted by $\bigwedge^k L_z$ and is defined as the $\Z$-span of the wedge products $v_1 \wedge \cdots \wedge v_k$ for all $v_1, \ldots, v_k \in L_z$. It is a lattice inside $\bigwedge^k V$. The inner product is given by
\begin{displaymath}
    \langle v_1 \wedge \cdots \wedge v_k, w_1 \wedge \cdots \wedge w_k \rangle = \det ( \langle v_i, w_j \rangle )_{1 \leq i,j \leq k}
\end{displaymath}
and extended linearly.

We have an isomorphism
\begin{displaymath}
    \bigwedge^{n-1} V \cong V,
\end{displaymath}
by sending $w \in \bigwedge^{n-1} V$ to $v \in V$ such that, for all $u \in V$,
\begin{displaymath}
    w \wedge u = \langle v, u \rangle.
\end{displaymath}
We make implicit use of the fact that $\bigwedge^n \R^n \cong \R$ and of an intermediary isomorphism with the dual space $V^\ast$. The isomorphism above is an isometry. 

Indeed, we can check that an orthonormal basis is sent to an orthonormal basis. Let $(e_1, \ldots, e_n)$ be the standard orthonormal basis of $V$. Then 
\begin{displaymath}
    (e_1 \wedge \cdots \wedge e_{n-1},\: e_1 \wedge \cdots \wedge e_{n-2} \wedge e_n,\: \ldots,\: e_2 \wedge \cdots \wedge e_n),
\end{displaymath}
is an orthonormal basis of $\bigwedge^{n-1} V$, formed by respectively removing each vector $e_i$ from the wedge product $e_1 \wedge \ldots \wedge e_n$. It is then easy to check that
\begin{displaymath}
    e_1 \wedge \ldots \wedge e_{n-1} \mapsto e_n, \; e_1 \wedge \ldots \wedge e_{n-2} \wedge e_n \mapsto -e_{n-1}, \; \ldots, \;  e_2 \wedge \ldots \wedge e_n \mapsto (-1)^{n-1} e_1.
\end{displaymath}

\begin{lemma} \label{lemma-isometric-lats}
    The lattice $\bigwedge^{n-1} L_z$ is isometric to the lattice $L_{\det(z) \cdot z^{-T}}$.
\end{lemma}
\begin{proof}
	We use the isomorphism $\bigwedge^{n-1} V \cong V$ described in the paragraphs above. 
	The wedge product has the property that $v_1 z \wedge \ldots \wedge v_n z = \det(z) \cdot v_1 \wedge \ldots \wedge v_n$ for $n$ row vectors $(v_i)$. This allows us to check that, under the given isomorphism,
	\begin{displaymath}
		e_1 z \wedge \ldots \wedge e_{n-1} z \mapsto \det(z) \cdot e_n z^{-T},
	\end{displaymath}
	and analogously for the other basis vectors above.
\end{proof}

\subsection{Successive minima}

Throughout this paper, we consider successive minima of lattices $L_z$ with respect to the unit ball $B^1 \subset V$ given by the standard inner product. When considering the exterior products of these lattices, successive minima are defined with respect to the compounds of the unit ball, as in the work of Mahler \cite{mahler} (refer also to \cite{evertse}, Section 3, for a modern treatment).

More precisely, the \emph{$k$-th compound} of $B^1$, denoted here by $B^k$, is defined as the convex hull of the points $x_1 \wedge \cdots \wedge x_k$, for all $x_1, \ldots, x_k \in B^1$. Mahler notes that $B^k$ is a bounded, convex body in $\bigwedge^k \R^n$, though generally not a sphere (see Section 4 in \cite{mahler}). Nevertheless, since $B^k$ is bounded and $0$ is an inner point of $B^k$, there are constants $c_{k,n}, C_{k,n} > 0$ such that
\begin{displaymath}
	B(n,k,c_{k,n}) \subset B^k \subset B(n,k,C_{k,n}),
\end{displaymath}
where $B(n,k,r)$ is the ball of radius $r$ inside $\bigwedge^k \R^n$. As such, the length $l$ of the shortest non-zero vector in $\bigwedge^k L_z$ can be approximated as 
\begin{displaymath}
	l \asymp_{n,k} \mu_1,
\end{displaymath}
where $\mu_1$ is the first successive minimum of $\bigwedge^k L_z$ with respect to $B^k$.

A theorem of Mahler (Theorem 3 in \cite{mahler}; Theorem 3.2 in \cite{evertse}) relates the successive minima of a lattice to those of its exterior powers. We state here a special case, relevant in this paper. 
	
\begin{lemma} \label{lemma-Mahler}
	Let $L$ be a lattice in $\mathbb{R}^n$ and let $\lambda_1 \leq \ldots \leq \lambda_n$ be its successive minima with respect to the unit ball $B^1$. Let $\mu_1$ be the first successive minimum of the lattice $\bigwedge^k L$ with respect to $B^k$. Then
	\begin{displaymath}
		\mu_1 \asymp_{n,k} \lambda_1 \cdots \lambda_k.
	\end{displaymath}
\end{lemma}

As explained above, this lemma implies that, if $l$ is the length of the shortest non-zero vector in $\bigwedge^k L$, then
\begin{displaymath}
	l \asymp_{n,k} \lambda_1 \cdots \lambda_k.
\end{displaymath}
We use this relation in Section \ref{sec-fricke-red}.

We also recall here a classical theorem of Minkowski (see \cite[Theorem VIII.1]{cassels}), stating that
\begin{equation} \label{eq-minkowski}
	d(L) \ll_n \lambda_1 \cdots \lambda_n \ll_n d(L),
\end{equation}
where $d(L)$ is the determinant of the lattice, e.g.\ $d(L_z) = \det(z)$. In particular, for a lattice of determinant $1$, called a \emph{unimodular lattice}, we have
\begin{equation} \label{eq-lambda-1}
	\lambda_1 \ll_n 1,
\end{equation}
using the inequalities $\lambda_1 \leq \lambda_i$, for all $i$.

The detailed study of successive minima of $L_z$ is crucial in this paper due to the following well-known lemma (see e.g. \cite[Lemma 1]{BHM}), which we apply when counting integral matrices, as explained at the end of Section \ref{sec-amp}.

\begin{lemma} \label{lemma-lattice-count}
	Let $L \subset \R^n$ be a lattice and let $\lambda_1 \leq \ldots \leq \lambda_n$ be its successive minima with respect to the unit ball. Let $B \subset \R^n$ be a ball of radius $R$ and arbitrary centre. We have the inequality
	\begin{displaymath}
		| L \cap B | \ll_{n} 1 + \frac{R}{\lambda_1} + \frac{R^2}{\lambda_1 \lambda_2} + \cdots + \frac{R^n}{\lambda_1 \cdots \lambda_n}.
	\end{displaymath}
\end{lemma}

\subsection{Iwasawa coordinates and reduction theory} \label{sec-Iwasawa}

Let $\uhp = \uhp_n$ be the generalised upper half plane, that is
\begin{displaymath}
	\uhp = \text{GL}_n(\R)/(\text{O}(n)\cdot \R^\times) \cong \SL_n(\R) / \text{SO}(n).
\end{displaymath}
In particular, the statement $z \in \uhp$ is taken to imply $z \in \SL_n(\R)$.
	
By the Iwasawa decomposition (see Section 1.2 in \cite{goldfeld}), we can take elements in $\uhp$ to be of the form $z = n(x) \cdot a(y)$, where $n(x) = (x_{ij})_{1 \leq i, j \leq n} \in \SL_n(\R)$ is upper triangular unipotent, meaning that it satisfies
\begin{displaymath}
	x_{ij} = \begin{cases}
		0, & j < i; \\
		1, & i = j;
	\end{cases}
\end{displaymath}
and $a(y)$ is diagonal, parametrised as
\begin{displaymath}
	a(y) = \diag(d_1, \ldots, d_n) = \diag(d y_1 \cdots y_{n-1},\; \ldots,\; d y_1 y_2,\; d y_1,\; d),
\end{displaymath}
where $d, y_1, \ldots, y_{n-1} \in \R_{>0}$ such that
\begin{displaymath}
	\det a(y) = d^{n} y_1^{n-1} y_2^{n-2} \cdots y_{n-1} = 1.
\end{displaymath}

\begin{remark} \label{rem-iwasawa-calculations}
	When considering the lattice that $z$ defines, the basis given by the Iwasawa decomposition will be particularly useful for computations.
	In particular, we note that the ``last'' basis element, or row, has the form 
	\begin{displaymath}
		e_n \cdot z = e_n \cdot n(x) a(y) = (0, 0, \ldots, 0, d)
	\end{displaymath}
	and the penultimate one is
	\begin{displaymath}
		e_{n-1} \cdot z = (0, \ldots, 0, dy_1, d x_{n-1, n}).
	\end{displaymath}
	Thus, $\norm{e_n}_z = d$ and, since there is only one $2\times 2$ minor with non-zero determinant in the matrix formed by the last two rows of $z$,
	\begin{displaymath}
		\norm{e_{n-1} \wedge e_n}_z = d^2 y_1
	\end{displaymath}
	Observe here also the useful computation
	\begin{displaymath}
		\norm{e_2 \wedge \cdots \wedge e_n}_z = \norm{e_1}_{z^{-T}} = (d y_1 \cdots y_{n-1})^{-1}.
	\end{displaymath}
\end{remark}

Define the Siegel set $\mathfrak{S}$ to be the set of all $z = n(x) a(y) \in \SL_n(\R)$ such that
\begin{displaymath}
	|x_{ij}| \leq \frac12
\end{displaymath}
for all $i < j$ and
\begin{displaymath}
	y_i \geq \frac{\sqrt{3}}{2},
\end{displaymath}
for all $i$, using the Iwasawa coordinates defined above.\footnote{In this paper, a Siegel set more generally refers to a set of points with $x$-coordinates bounded from above in absolute value and $y$-coordinates bounded from below.}
Reduction theory (see \cite[Theorem I.1.4]{borel} or \cite[Proposition 1.3.2]{goldfeld}) shows that
\begin{displaymath}
	\SL_n(\R) = \SL_n(\Z) \cdot \mathfrak{S}.
\end{displaymath}
If $z \in \mathfrak{S}$, we say that $(e_1 z, \ldots, e_n z)$ is a \emph{reduced basis} for $L_z$. We also remark that reduction theory allows us to pick $e_n z$ to be any vector of shortest length in $L_z$ (this is, indeed, part of the reduction algorithm).

\begin{remark} \label{rem-reduction-n-1}
	It is useful in later sections to note an embedding of $\SL_{n-1}(\R)$ into $\SL_n(\R)$ and the connection between the two systems of Iwasawa coordinates. More precisely, we can write $z = n(x) a(y) \in \mathbb{H}$ as 
	\begin{displaymath}
		z = \begin{pmatrix}
			dy_1 \cdot w & \ast \\
			0 & d
		\end{pmatrix},
	\end{displaymath}
	where $w \in \mathbb{H}_{n-1}$ is a matrix in $\GL_{n-1}(\R)$. Though not normalised, we can use a variant of the Iwasawa coordinates (it is the one used in Definition 1.2.3 in \cite{goldfeld}) to write $w = n(x') \cdot a(y')$, where
	\begin{displaymath}
		a(y') = \diag(y_2 \cdots y_{n-1}, \ldots, y_2, 1).
	\end{displaymath}
	Multiplication of $z$ by parabolic matrices 
	\begin{displaymath}
		g = \begin{pmatrix}
			h & 0 \\
			0 & 1
		\end{pmatrix} \in \SL_n(\Z)
	\end{displaymath}
	with $h \in \SL_{n-1}(\Z)$, acts on $w$ by sending it to $h \cdot w$ and otherwise leaves the last row of $z$ invariant. Reduction theory in degree $n-1$ now implies that there is a parabolic block matrix $g \in \SL_n(\Z)$ as above so that $g \cdot z = n(x) \cdot a(y)$ with $y_i \geq \sqrt{3}/2$ for $i = 2, \ldots, n-1$.
\end{remark}

We recall also another standard lemma, which informally says that a reduced basis behaves similarly to an orthogonal basis.

\begin{lemma} \label{lemma-red-basis-orthog}
	Let $(v_1, \ldots, v_n)$ be a reduced basis of a lattice $L$. Let $v \in L$ and write $v = \sum_{i=1}^{n} a_i v_i$ with $a_i \in \Z$. Then $a_i \ll_n \norm{v} / \norm{v_i}$.
\end{lemma}
\begin{proof}
	See Lemma 1 in \cite{venk}.
\end{proof}

Finally, if $\Omega \subset \mathbb{H}$ is a compact set (in particular, it projects to a compact set in the space of lattices $\SL_n(\Z) \backslash \mathbb{H}$) and $z \in \Omega$, then $\lambda_1 \gg_\Omega 1$ by Mahler's criterion \cite[Corollary I.1.9]{borel}. The other successive minima must then also be bounded from below, so $\lambda_i \gg 1$. By \eqref{eq-minkowski}, we have that 
\begin{displaymath}
	1 \ll \lambda_2^{n-1} \leq \lambda_2 \cdots \lambda_n \ll 1/\lambda_1 \ll 1
\end{displaymath}
since $z$ has determinant $1$. Thus, $\lambda_2 \asymp 1$ and inductively we find $\lambda_i \asymp_\Omega 1$ for all $i$. We may say $L_z$ is an $\Omega$-balanced lattice.

For any $z \in \mathbb{H}$ we say that $z$ \emph{reduces} to $\Omega$ if there is $w \in \Omega$ such that $L_z = L_w$, in other words if there is $\gamma \in \SL_n(\Z)$ such that $z = \gamma w$. The discussion in the paragraph above proves the following lemma.

\begin{lemma} \label{lemma-red-to-cpt}
	Suppose that $z \in \mathbb{H}$ reduces to a compact set $\Omega$ and let $\lambda_1, \ldots, \lambda_n$ be the successive minima of $L_z$. Then $\lambda_i \asymp_{\Omega} 1$ for all $i \in \{1, \ldots, n\}$, where the implicit constant depends only on $\Omega$.
\end{lemma}

For some parameter $\alpha$, we call a lattice $L$ $\alpha$-\emph{balanced} if $\lambda_1 \asymp_\alpha d(L)^{1/n}$, where $\lambda_1 \leq \ldots \leq \lambda_n$ are the successive minima of $L$.
The lemma above shows that $L_z$ is $\Omega$-balanced if $z$ reduces to $\Omega$. 
Since $\Omega$ will always be fixed, we often call $\Omega$-balanced lattices simply \emph{balanced}.

\begin{lemma} \label{lemma-balancedness}
	A lattice $L_z$ is $\alpha$-balanced if and only if $\lambda_i \asymp_\alpha d(L)^{1/n}$ for all $i = 1, \ldots, n$.
	Moreover, this property is invariant under scaling and taking duals.
\end{lemma}
\begin{proof}
	This follows from Minkowski's theorem \eqref{eq-minkowski}, generalising Lemma \ref{lemma-red-to-cpt}.
	From Lemma \ref{lemma-isometric-lats} on the dual lattice and Lemma \ref{lemma-Mahler} on the successive minima of exterior products, we deduce that $L$ is balanced if and only if the dual $L^\ast$ is balanced. 
\end{proof}

\section{The amplified pre-trace formula} \label{sec-amp}
	
We follow the amplification scheme of Blomer and Maga \cite{BM4}, using their archimedean test function but giving also a version that simplifies the sum over Hecke eigenvalues by assuming a conjecture about their sizes. 

Let $G = \SL_n(\R)$, $K = \SO(n)$, $\Gamma = \Gamma_0(N)$, and let $\phi$ be the cuspidal Hecke-Maaß form of level $N$ that we wish to bound. 
Let $\mu = (\mu_1, \ldots, \mu_n)$ be the spectral parameters of $\phi$. 
We may embed $\phi$ into a basis of the space of Hecke-Maaß cusp forms for $\Gamma_0(N)$. 
More precisely, we have a spectral decomposition
\begin{displaymath}
	L^2(\Gamma_0(N) \backslash \mathbb{H}) = \int V_\varpi \, d\varpi = L^2_{\text{cusp}} \oplus L^2_{\text{Eis}},
\end{displaymath}
where every $V_\varpi$ is a one-dimensional space generated by an eigenform $\phi_\varpi$ of the algebra of invariant differential operators and the Hecke algebra. 
Let $\mu_\varpi$ be the spectral parameter of $\phi_\varpi$ and assume that $\phi = \phi_{\varpi_0}$. 
Note moreover that $L^2_{\text{cusp}}$ has a discrete decomposition.

Recall the Cartan decomposition $G = KAK$, where $A$ is the subgroup of diagonal matrices. 
The latter has a Lie algebra $\mathfrak{a}$, on which the Weyl group $W$ of $G$ acts. 
We define the Cartan projection $C(g) \in \mathfrak{a}/W$ of an element $g \in G$ via the Cartan decomposition $g = k_1 \exp(C(g)) k_2$, where $k_1, k_2 \in K$. 
Now pick a $W$-invariant norm $\norm{\cdot}$ on $\mathfrak{a}$. 
We note that, if $\norm{C(g)} \ll 1$, then by exponentiating we have
\begin{displaymath}
	g = k + O(1),
\end{displaymath}
where $k \in K$ and $O(1)$ stands for a matrix whose norm (by equivalence, any norm) is $O(1)$.

\subsection{The Hecke algebra and Hecke eigenvalues} \label{sec-hecke-alg}

We now briefly review some aspects of the structure of the unramified Hecke algebra. Let $p$ be a prime not dividing $N$ and $\textbf{a} = (a_1, \ldots, a_n) \in \Z^n$. The double coset
\begin{displaymath}
\Gamma \diag(p^{a_1}, \ldots, p^{a_n}) \Gamma = \bigcup_j \Gamma \alpha_j
\end{displaymath}
defines a Hecke operator
\begin{displaymath}
T_{\textbf{a}}(p) (\psi)(z) = \sum_j \psi(\alpha_j \cdot z),
\end{displaymath}
where $\psi$ is any function on $\Gamma \backslash \mathbb{H}$. 
We define the standard Hecke operator as
\begin{displaymath}
	T(p) = T_{(1, 0, \ldots, 0)}(p).
\end{displaymath}
One computes that the adjoint of $T(p)$ is the operator $T'(p) = T_{(1, \ldots, 1, 0)}(p)$.
Let $\lambda(p, \phi_\varpi)$ be the eigenvalue of $\phi_\varpi$ under $T(p)$, so that $\overline{\lambda(p, \phi_\varpi)}$ is its eigenvalue under $T'(p)$.
By \cite[Lemma 4.4]{BM4} we have
\begin{displaymath}
T(p) \cdot T'(p) = a \cdot T_{(2, 1, \ldots, 1, 0)}(p) + b \cdot p^{n-1} \operatorname{id},
\end{displaymath}
where $a, b \ll 1$. 
Furthermore, if $p$ and $q$ are distinct primes not dividing $N$, then we have the multiplication rule on double cosets
\begin{equation} \label{eq-hecke-sets}
\Gamma \diag(p, 1, \ldots, 1) \Gamma \cdot \Gamma \diag(q, \ldots, q, 1) \Gamma = \Gamma \diag(pq, q, \ldots, q, 1) \Gamma
\end{equation}
in the Hecke algebra, corresponding to the composition $T(p) \cdot T'(q)$ (see \cite[Section 6]{BM4}).

Let now $L > 0$ be a parameter and $\mathcal{P}$ be the set of primes contained in $[L, 2L]$, not dividing $N$. Define 
\begin{displaymath}
	A_\varpi = \left| \sum_{p \in \mathcal{P}} \frac{\lambda(p, \varpi)}{p^{(n-1)/2}} \cdot x_p \right|^2,
\end{displaymath}
where $x_p = |\lambda(p, \varpi_0)|/\lambda(p, \varpi_0)$.

We use here the normalised eigenvalues $\lambda(p, \varpi)/p^{(n-1)/2}$ as defined in \cite[(9.3.5)]{goldfeld}.
Note that 
\begin{displaymath}
	A_{\varpi_0} = \left| \sum_p |\lambda(p, \varpi_0)/p^{(n-1)/2}| \right|^2.
\end{displaymath}
A lower bound for this quantity is given in Hypothesis \eqref{eq-hypothesis}. 
We now prove it follows from GRH.

\begin{lemma} \label{lemma-hypothesis}
	Let $\delta > 0$ be any positive constant and $N \gg_\delta 1$ be large enough. 
	Assuming the Grand Riemann Hypothesis, if $L > N^\delta$, then
	\begin{equation}
		\sum_{p \in \mathcal{P}} \frac{|\lambda(p, \varpi_0)|}{p^{(n-1)/2}} \gg_{\varepsilon} L^{3/4-\varepsilon}.
	\end{equation}
\end{lemma}
\begin{proof}
	The following are standard computations, and we refer to Sections 5.1, 5.3, 5.6, 5.7 in \cite{IK} for more details. 
	Let $\lambda(p) = \lambda(p, \varpi_0)/p^{(n-1)/2}$ and note that these give the coefficients of the $L$-function attached to $\phi$ or, equivalently, to the automorphic representation $\pi$ generated by $\phi$.
	Let $L_{\text{RS}}(s) = L(s, \pi \times \tilde{\pi})$ be the Rankin-Selberg $L$-function and define $\Lambda_{\text{RS}}(n)$ to be its coefficients, so that
	\begin{displaymath}
		\frac{L'_{\text{RS}}}{L_{\text{RS}}}(s) = \sum_{n = 1}^\infty \frac{\Lambda_{\text{RS}}(n)}{n^s}.
	\end{displaymath}
	Then we have $\Lambda_{\text{RS}}(p) = |\lambda(p)|^2 \log p$.
	
	The prime number theorem under GRH states that
	\begin{equation} \label{eq-pnt}
		\sum_{n \leq x} \Lambda_{\text{RS}}(n) = x + O_{\varepsilon, \mu}(x^{1/2 + \varepsilon} \cdot N^{\varepsilon}).
	\end{equation}
	For $y \leq \sqrt{x}$, we obtain that
	\begin{displaymath}
		\sum_{x\leq n \leq x + y} \Lambda_{\text{RS}}(n) \ll_{\varepsilon, \mu} x^{1/2 + \varepsilon} N^\varepsilon.
	\end{displaymath}
	
	Now we note that $\Lambda_{\text{RS}}(n) \geq 0$ for all $n$ by the definition of the Rankin-Selberg convolution. 
	It follows from the prime number theorem above by dropping all but one term that
	\begin{displaymath}
		\lambda(p)^2 \ll \Lambda_{\text{RS}}(p) \ll x^{1/2 + \varepsilon} N^\varepsilon 
	\end{displaymath}
	for $p \asymp x$.
	
	Let $x \gg N^\delta$ for some $\delta > 0$. 
	The bound above and \eqref{eq-pnt} imply that
	\begin{displaymath}
		x^{1 - \varepsilon} \ll \sum_{p \asymp x} |\lambda(p)|^2 \ll x^{1/4 + \varepsilon} \sum_{p \asymp x} |\lambda(p)|.
	\end{displaymath}
	This proves the claim.
\end{proof}

\begin{remark} \label{rem-ramanujan}
	It is expected that a stronger version of \eqref{eq-hypothesis} holds, that is, with exponent $1$ instead of $3/4$.
	To prove this we would require the Ramanujan-Petersson conjecture. 
	This would improve the saving in Theorem \ref{thm-main} by doubling the exponent.
\end{remark}

\subsection{Amplifiers}
Looking towards applying the amplified pre-trace formula, we choose the archimedean test function $f_\mu: C_c^\infty(K \backslash G / K) \longrightarrow \mathbb{C}$ defined in \cite[Section 3]{BM4}. 
It has compact support, independent of $\mu$, and is bounded $f_\mu \ll_{\mu,n} 1$ in terms of $\mu$, where the dependence on $\mu$ is continuous.\footnote{In fact, there is an explicit bound for the function $f_\mu$ that would recover the baseline bound in the spectral aspect in our main theorem. For our purposes, we may simply bound $f_\mu$ by a constant depending on $\mu$ and $n$, but independent of the level.}

We shall consider matrices $g \in G$ in the support of $f_\mu$.
Such an assumption implies that $g = k + O_n(1)$ for some $k \in K$, where $O_n(1)$ refers to a matrix with entries bounded by some constant depending on $n$.
Now, let $\mathcal{M}_n(\Z,N)$ be the set of integral matrices with last row congruent to $(0, \ldots, 0,*)$ modulo $N$. For $(m, N) = 1$, let
\begin{displaymath}
	H(z, m, N) := \{ \gamma \in \mathcal{M}_n(\Z,N) \mid \det \gamma = m,\, z^{-1} \gamma z = O( m^{1/n}) \},
\end{displaymath}
where $z^{-1} \gamma z = O( m^{1/n})$ means that $m^{-1/n} z^{-1} \gamma z$ lies in the support of $f_\mu$, as above.
The implicit constant thus depends on $n$ and the support of $f_\mu$, dependence which we suppress throughout the paper.

\begin{proposition} \label{prop-amp}
	Let $\phi$ be a Hecke-Maaß form for $\Gamma_0(N) \leq \SL_n(\R)$ with spectral parameter $\mu$, let $L \gg N^\delta$ for some $\delta > 0$ be a parameter and let $\mathcal{P}$ be the set of primes in $[L, 2L]$, not dividing $N$. Then, assuming Hypothesis \eqref{eq-hypothesis}, we have the bound
	\begin{displaymath}
		L^{3/2-\varepsilon} |\phi(z)|^2 \ll_{\mu, \varepsilon} |\mathcal{P}| \cdot |H(z, 1, N)| + \frac{1}{L^{n-1}}  \sum_{p, q \in \mathcal{P}} |H(z, p \cdot q^{n-1}, N)|.
	\end{displaymath}
\end{proposition}
\begin{proof}
	In preparation for the proof, recall that, by the work of Blomer and Maga, the spherical transform $\tilde{f}_\mu$ of $f_\mu$ satisfies $\tilde{f}_\mu(\mu) \geq 1$ and is non-negative on all possible spectral parameters occurring in the decomposition of $L^2(\Gamma_0(N) \backslash \mathbb{H})$. 
	Additionally, when writing $f_\mu(g)$ for $g \in \GL_n(\R)$, where $\det(g) > 0$, we mean $f_\mu(g / \det(g)^{1/n})$, thus extending the domain of $f_\mu$ by postulating its invariance under scalars.
	
	Now consider
	\begin{displaymath}
		\int A_{\varpi} \cdot \tilde{f}_\mu(\mu_\varpi) \phi_\varpi(z) \overline{\phi_\varpi(w)}\, d\varpi,
	\end{displaymath}
	expand every $A_\varpi$ and group terms into expressions of the form
	\begin{displaymath}
		\frac{1}{(pq)^{(n-1)/2}} \int \lambda(p, \varpi) x_p \cdot \overline{\lambda(q, \varpi) x_q} \cdot \tilde{f}_\mu(\mu_\varpi) \phi_\varpi(z) \overline{\phi_\varpi(w)}\, d\varpi,
	\end{displaymath}
	which is equal to 
	\begin{displaymath}
	S_{p,q} = \frac{x_p \overline{x_q}}{(pq)^{(n-1)/2}} \cdot T(p) T'(q) \cdot  \int \tilde{f}_\mu(\mu_\varpi) \phi_\varpi(z) \overline{\phi_\varpi(w)}\, d\varpi,
	\end{displaymath}
	where the Hecke operators act in the variable $z$. 
	We apply the pre-trace formula to obtain the geometric side
	\begin{displaymath}
		S_{p,q} = \frac{x_p \overline{x_q}}{(pq)^{(n-1)/2}}  \cdot T(p) T'(q) \sum_{\gamma \in \Gamma} f_\mu(z^{-1} \gamma w),
	\end{displaymath}
	where again we write $\Gamma_0(N) = \Gamma$ for brevity.
	Note that for any double coset $\Gamma g \Gamma$, the corresponding Hecke operator $T_g$ acts on the variable $z$ by
	\begin{displaymath}
		T_g \sum_{\gamma \in \Gamma} f_\mu(z^{-1}\gamma w) = \sum_{\gamma \in \Gamma g \Gamma} f_\mu(z^{-1} \gamma w), 
	\end{displaymath}
	by definition and sum unfolding. Moreover, using the compact support of $f_\mu$, we can bound the right-hand side by
	\begin{displaymath}
		\sum_{\gamma \in \Gamma g \Gamma} f_\mu(z^{-1} \gamma w) \ll_\mu | \{ \gamma \in \Gamma g \Gamma \mid z^{-1} \gamma w = \det(\gamma)^{1/n} (k + O(1)),\, k \in K \} |
	\end{displaymath}
	using the triangle inequality. Since $K$ is compact, we can simplify $k + O(1)$ to $O(1)$, where the implicit constant depends on $n$.
	
	We now write the compositions $T(p) \cdot T'(q)$ as linear combinations of Hecke operators $T_g$. Let $z = w$ and assume that $p \neq q$. Recalling that $T(p) \cdot T'(q)$ is the Hecke operator corresponding to 
	\begin{displaymath}
	\Gamma \diag(pq, q, \ldots, q, 1) \Gamma,
	\end{displaymath}
	and that $x_p \ll 1$ for all $p \in \mathcal{P}$, we bound
	\begin{displaymath}
	S_{p,q} \ll_\mu \frac{1}{L^{n-1}} \cdot | H(z, p q^{n-1}, N) |.
	\end{displaymath}
	Note that we have made this upper bound larger by forgetting the structure of the double coset and simply retaining the information about the determinant, which is an invariant of the double coset.
	Analogously we obtain
	\begin{displaymath}
		S_{p,p} \ll_\mu \frac{1}{L^{n-1}} \cdot | H(z, p^{n}, N) | + | H(z, 1, N) |.
	\end{displaymath}
	
	We now put together the bounds above and observe that non-negativity of $\tilde{f}_\mu$ and of $A_\varpi$ gives
	\begin{displaymath}
	A_{\varpi_0} |\phi(z)|^2 \leq \int A_{\varpi} \cdot \tilde{f}_\mu(\mu_\varpi) |\phi_\varpi(z)|^2 d\varpi.
	\end{displaymath}
	Finally, we get a lower bound on $A_{\varpi_0}$ by Hypothesis \eqref{eq-hypothesis}.
\end{proof}

For unconditional bounds, one may work with the amplifier given in \cite[(6.2)]{BM4}. 
It uses Hecke operators attached to higher powers of primes for providing an alternative to Hypothesis \eqref{eq-hypothesis}. 
In fact, we give the slightly more precise version of this amplifier by including information on the determinantal divisors. 
Recall that the \emph{$j$-th determinantal divisor} $\Delta_j(\gamma)$ of an integral matrix $\gamma$ is equal to the greatest common divisor of all $j \times j$ minors.

\begin{proposition} \label{prop-amp-uncond}
	With the same notation as in Proposition \ref{prop-amp}, we have the unconditional bound
	\begin{displaymath}
	L^{2-\varepsilon} |\phi(z)|^2 \ll_{\mu, \varepsilon} |\mathcal{P}| \cdot |H(z, 1, N)| + \sum_{\nu = 1}^n \frac{1}{L^{(n-1)\nu}} \sum_{p, q \in \mathcal{P}} |\overline{H}(z, p^\nu, q^{(n-1)\nu}, N)|,
	\end{displaymath}
	where $\overline{H}(z, p^\nu, q^{(n-1)\nu}, N)$ consists of matrices $\gamma \in H(z, p^\nu q^{(n-1)\nu}, N)$ satisfying the additional conditions
	\begin{displaymath}
		\Delta_j(\gamma) = (q^\nu)^{j-1},
	\end{displaymath}
	for all $1 \leq j \leq n-1$.
\end{proposition}
\begin{remark}
	Blomer and Maga only preserve the condition on $\Delta_1$ and $\Delta_2$ (see their definition of $S(m,l)$). 
	These and the additional ones in the proposition above follow directly using the crucial property of the determinantal divisors, namely their invariance under right or left multiplication by elements of $\SL_n(\Z)$ (see e.g.\ \cite[Thm. II.8]{newman}). 
	Except for the proof of Proposition \ref{prop-degenerate-count}, these conditions are not used, and we mostly consider the larger set $H(z,m,N)$ for simplicity of notation.
\end{remark}

\section{Higher rank Atkin-Lehner operators} \label{sec-fricke}

In this section we consider possible generalisations of Atkin-Lehner operators to the spaces $X_n(N)$ for $n > 2$.
We consider this to be of independent interest and therefore do a thorough investigation of all cases, regardless of the restrictions imposed in the rest of this paper.
In fact, the results in this section motivate these restrictions, as one of the main conclusions is the uniqueness of the generalised Fricke involution among the potential symmetries of $X_n(N)$ considered here for $n > 2$.

\subsection{The normaliser of the Hecke congruence subgroup} \label{sec-normaliser}

In the theory of automorphic forms on $\SL_2(\R)$, an Atkin-Lehner operator $S$ is an involution on space of left-$\Gamma_0(N)$ invariant functions.
It is obtained by setting $Sf (z) = f (gz)$ for all $z \in \mathbb{H}$, where $g$ lies in the normaliser of $\Gamma_0 (N)$ inside $\SL_2(\R)$. 
This is a natural method of producing automorphisms, since the invariance of $f (z)$ under a group $\Gamma$ is equivalent to the invariance of $f (gz)$ under $g^{-1} \Gamma g$. 
The normaliser has been computed by Atkin and Lehner in \cite{AL} and an example of a non-trivial normalising element is
\begin{displaymath}
	g = \begin{pmatrix}
		& - 1\\
		N & 
	\end{pmatrix},
\end{displaymath}
which induces the so-called \emph{Fricke involution}.
In fact, the normaliser gives all automorphism of the modular curve $X_2(N)$, in more standard notation $X_0(N)$, for all $N$ up to finitely many exceptions (see \cite{KM}).

Thus, searching for symmetries of automorphic forms in higher rank should involve computing the normalisers of $\Gamma_0 (N) \leq \SL_n(\R)$ for $n > 2$.
Unfortunately, this method can only produce the identity operator, since we prove below that these normalisers, in contrast to the case $n = 2$, are trivial. 
In the following we denote by $\GL_n^+ (\Q)$ the subgroup of invertible matrices with positive determinant.

\begin{theorem} \label{thm-normaliser}
	For $n > 2$, the normaliser of $\Gamma_0 (N)$ inside $\GL_n^+(\Q)$ is trivial, that is, equal to $\mathbb{Q}_{> 0} \cdot \Gamma_0(N)$.
\end{theorem}

For simplicity and clarity of the argument, since we work with some explicit coordinates, we prove the theorem in the case of $n = 3$.
The way to generalise the proof should be apparent to the reader.

Consider the \emph{left} action of $G := \GL_3^+(\Q)$ on full $\Z$-lattices in $\R^3$ (using column vectors).\footnote{As opposed to the rest of the present paper, in this independent section we let $G$ act from the left on vectors. This allows for some simplifications of the arguments. In fact, from the point of view of lattices, this is the more natural setting for $\Gamma_0(N)$. For instance, when $N$ is prime, it is easier to see that $\SL_n(\R) / \Gamma_0(N)$ parametrises pairs of unimodular lattices together with a sublattice of index $N$. On the other hand, in the theory of automorphic forms, the dual picture is more standard.} 
Let $L_1 = \langle e_1, e_2, e_3 \rangle$ be the standard lattice for a basis $(e_1, e_2, e_3)$ of $\R^3$ and consider $\mathcal{L}= G \cdot L_1$, the orbit of $L_1$ under the action of $G$.

Note that the stabiliser of $L_1$ under this action is the group $\SL_3(\Z)$. 
More generally, for $M \in \N$, let $L_M = \langle e_1,e_2, M e_3 \rangle$, or in other words,
\begin{displaymath} 
	L_M = \begin{pmatrix}
		1 &  & \\
		& 1 & \\
		&  & M
	\end{pmatrix}
	\cdot L_1.
\end{displaymath}
If we let $A_M = \diag(1,1,M)$, then the stabiliser of $L_M$ is
\begin{displaymath}
	\operatorname{Stab}(L_M) = A_M \operatorname{Stab}(L_1) A_M^{- 1} = \left\{ 
	\begin{pmatrix}
		a_{11} & a_{12} & \frac{a_{13}}{M}\\
		a_{21} & a_{22} & \frac{a_{23}}{M}\\
		Ma_{31} & Ma_{32} & a_{33}
	\end{pmatrix}
	: (a_{ij}) \in \SL_3(\mathbb{Z}) \right\} .
\end{displaymath}
It follows that $\operatorname{Stab}(L_1) \cap \operatorname{Stab}(L_M) = \Gamma_0(M)$.
Since $\Gamma_0(N) \subset \Gamma_0(M)$ for all $M \mid N$, we also have that
\begin{displaymath}
	\bigcap_{M \mid N} \stab(L_M) = \Gamma_0(N).
\end{displaymath}
The following lemma provides a converse for this observation.

\begin{lemma} \label{lemma-fixed-lats}
	The set of lattices fixed by $\Gamma_0 (N)$ is
	\begin{displaymath}
		\bigcup_{M \mid N} \{ q L_M : q \in \Q_{> 0} \}.
	\end{displaymath}
\end{lemma}

\begin{proof}
	Let $L = g \cdot L_1 \in \mathcal{L}$, where $g \in \GL_3^+(\Q)$, and assume that $\Gamma_0(N)$ fixes $L$. Then $g^{- 1}\Gamma_0(N) g$ fixes $L_1$, so we must have $g^{- 1} \Gamma_0 (N) g \subset \SL_3(\Z)$.
	
	Scaling $g$ by a positive rational number, we may assume that $g \in \mathcal{M}_{3 \times 3}(\mathbb{Z})$.
	Let then $H$ be the Hermite normal form of $g$, so that
	\begin{displaymath}
		H = gU, 
	\end{displaymath}
	with $U \in \SL_3(\Z)$ and $H$ lower triangular. 
	We have $H L_1 = gUL_1 = gL_1 = L.$ 
	So we may further assume that $g = H$ and is thus lower triangular. More explicitly, write
	\begin{displaymath}
		H = \begin{pmatrix}
			\alpha_1 & 0 & 0\\
			\beta_1 & \beta_2 & 0\\
			\gamma_1 & \gamma_2 & \gamma_3
		\end{pmatrix}
		\in \mathcal{M}_{3 \times 3}(\Z). 
	\end{displaymath}
	
	We test the inclusion $H^{- 1} \xi H \in \SL_3(\Z)$ with various matrices $\xi \in \Gamma_0(N)$. Observe that
	\begin{align*}
		H^{- 1} \begin{pmatrix}
			1 & 1 & \\
			& 1 & \\
			&  & 1
		\end{pmatrix} H \in \SL_3 (\mathbb{Z}) \quad
		& \text{implies that }
		& \frac{\beta_1}{\alpha_1}, \frac{\beta_2}{\alpha_1},
		\frac{\beta_1 \gamma_2 - \gamma_1 \beta_2}{\alpha_1 \gamma_3} \in
		\mathbb{Z};\\
		H^{- 1} \begin{pmatrix}
			1 &  & 1\\
			& 1 & \\
			&  & 1
		\end{pmatrix} H \in \SL_3 (\mathbb{Z}) \quad
		& \text{implies that } 
		& \frac{\gamma_1}{\alpha_1}, \frac{\gamma_2}{\alpha_1},
		\frac{\gamma_3}{\alpha_1} \in \mathbb{Z};\\
		H^{- 1} \begin{pmatrix}
			1 &  & \\
			1 & 1 & \\
			&  & 1
		\end{pmatrix} H \in \SL_3 (\mathbb{Z})  \quad
		& \text{implies that } 
		& \frac{\alpha_1}{\beta_2}, \frac{\alpha_1}{\beta_2} \cdot
		\frac{\gamma_2}{\gamma_3} \in \mathbb{Z};\\
		H^{- 1} \begin{pmatrix}
			1 &  & \\
			& 1 & \\
			N &  & 1
		\end{pmatrix} H \in \SL_3 (\mathbb{Z}) \quad
		& \text{implies that } 
		& N \frac{\alpha_1}{\gamma_3} \in \mathbb{Z}.
	\end{align*}
	Since $\frac{\beta_2}{\alpha_1}, \frac{\alpha_1}{\beta_2} \in \mathbb{Z}$,
	we must have $\frac{\beta_2}{\alpha_1} = \pm 1$. Since
	$\frac{\gamma_3}{\alpha_1}, N \frac{\alpha_1}{\gamma_3} \in \mathbb{Z}$, we
	must have $\frac{\gamma_3}{\alpha_1} = \pm M$, where $M \mid N$. Using the rest of the
	findings above, we may do column manipulations and obtain
	\begin{displaymath}
		H = \alpha_1 \begin{pmatrix}
			1 & 0 & 0\\
			\frac{\beta_1}{\alpha_1} & \frac{\beta_2}{\alpha_1} & 0\\
			\frac{\gamma_1}{\alpha_1} & \frac{\gamma_2}{\alpha_1} &
			\frac{\gamma_3}{\alpha_1}
		\end{pmatrix} = \alpha_1 \begin{pmatrix}
			1 &  & \\
			& 1 & \\
			&  & M
		\end{pmatrix} U',
	\end{displaymath}
	with $U' \in \SL_3 (\mathbb{Z})$. 
	Thus $L = H L_1 = L_M$ up to $\mathbb{Q}_{> 0}$ scalars.
\end{proof}

\begin{proof}[Proof of Theorem \ref{thm-normaliser}]
	Let $g \in \GL_3^+ (\mathbb{Q})$ such that $g^{- 1} \Gamma_0 (N) g = \Gamma_0 (N)$. 
	Since $\Gamma_0 (N)$ fixes the lattices $L_M$ for all divisors $M$ of $N$, we find that $\Gamma_0 (N)$ must also fix the lattices $g L_M$ for $M \mid N$. 
	By the previous lemma, for each divisor $M$ of $N$ there is a rational number $q_M$ and a divisor $f (M) \mid N$ such that
	\begin{displaymath}
		gL_M = q_M L_{f (M)}
	\end{displaymath}
	for all $M \mid N$.
	
	By the definition of $L_M$ and using the fact that $\stab  (L_1) =
	\SL_3 (\mathbb{Z})$, we can deduce that
	\begin{equation} \label{eq-mat}
		q_M^{- 1} \begin{pmatrix}
			1 &  & \\
			& 1 & \\
			&  & f (M)^{- 1}
		\end{pmatrix} \cdot g \cdot \begin{pmatrix}
			1 &  & \\
			& 1 & \\
			&  & M
		\end{pmatrix} \in \SL_3 (\mathbb{Z}),
	\end{equation}
	for all $M|N$.
	
	Rescaling $g$ by $q_1 \in \mathbb{Q}$ we may assume that $q_1 = 1$. 
	Taking $M = 1$ in \eqref{eq-mat} and applying determinants, we deduce that $\det(g) = f (1)$. 
	Applying determinants to all other equations, we find that
	\begin{displaymath}
		q_M^3 = \frac{f (1) M}{f (M)}. 
	\end{displaymath}
	In particular, for $M = N$, we have $q_N^3 f (N) = N f (1)$. 
	Since $f (N) \mid N$, we must have $q_N \in \mathbb{Z}$.
	
	Let us make \eqref{eq-mat} more explicit. Taking $M = 1$, we have
	\begin{displaymath}
		g = \begin{pmatrix}
			\ast & \ast & \ast\\
			\ast & \ast & \ast\\
			f (1) \ast & f (1) \ast & f (1) \ast
		\end{pmatrix},  
	\end{displaymath}
	where $\ast$ denotes unknown \emph{integers}. 
	In particular, the last column of $g$ is integral. 
	If we now take $M = N$, we have
	\begin{displaymath}
		g = \begin{pmatrix}
			q_N \ast & q_N \ast & \ast\\
			q_N \ast & q_N \ast & \ast\\
			q_N f (N) \ast & q_N f (N) \ast & \ast
		\end{pmatrix}.
	\end{displaymath}
	Using the properties of the determinant and that $\ast$ denotes integers, we
	deduce that $q_N^2 | \det (g) = f (1)$.
	
	Let $f (1) = q_N^2 k$ for some $k \in \mathbb{Z}$. 
	Now the last row of $g$ is divisible by $q_N^2 k$ and the first two columns are divisible by $q_N$.
	By the same method we infer that $q_N  k \cdot q_N \cdot q_N = q_N^3 k$ divides $\det (g) = f (1) = q_N^2 k$. 
	Therefore $q_N = 1$, which implies that $f (N) = N f (1)$. 
	Since $f (N) \mid N$, it follows that $f (1) = 1$ and $f (N) = N$. 
	Putting everything together, it follows that $g \in \Gamma_0 (N)$.
\end{proof}

\begin{remark}
	The case $n > 3$ can be done similarly. 
	In essence, what makes the case $n > 2$ differ from $n = 2$ is the imbalance between the number of columns with divisibility conditions and the number of rows with such conditions. 
	This leads to the different exponents of $q_N$ in the proof and ultimately to the triviality of the solutions to our equations.
\end{remark}

Theorem \ref{thm-normaliser-R} on the normaliser of $\Gamma_0(N)$ in the real group $\PGL_n(\R)$ now follows as a corollary to Theorem \ref{thm-normaliser}.
\begin{proof}[Proof of Theorem \ref{thm-normaliser-R}]
	We use the results of \cite{borel-dens}, which imply that the normaliser of $\Gamma_0(N)$, being commensurable with the arithmetic group $\PGL_n(\Z)$, lies in $\PGL_n(\Q)$.
\end{proof}

\subsection{A different perspective}

We have seen in the last section that $n = 2$ is singular in the sequence of families $\Gamma_0(N) \leq \SL_n(\Z)$ of congruence subgroups. 
To arrive at a general definition of Atkin-Lehner operators, it is useful to note another way in
which the group $\PGL_2 (\mathbb{R})$ is distinguished, as described below.

An important automorphism of matrices in $\SL_n (\mathbb{R})$ is the map $g \mapsto g^{- T}$, sending a matrix to its inverse transpose. 
As already noted in the present paper, this map sends a lattice $L_g$ to its dual, but is also used to define the dual form of an automorphic form for $\SL_n (\mathbb{Z})$ (see section 9.2 in \cite{goldfeld}) or also the contragredient representation of a $\GL(n)$ automorphic representation. 

In $\PGL(2)$, dual forms are not commonly mentioned because this automorphism is, in fact, inner in this case. 
Indeed, if we take 
\begin{displaymath}
	w =
	\begin{pmatrix}
		& - 1\\
		1 & 
	\end{pmatrix}
\end{displaymath}
to be the non-trivial Weyl element, then we easily compute that
\begin{equation}
	wg^{- T} w^{- 1} = - \frac{1}{\det (g)} g. \label{eq-n=2}
\end{equation}
In particular, the map $z \mapsto z^{- T}$ induces the identity on $\PGL_2(\Z) \backslash \PGL_2(\R) / \operatorname{PO}(2)$.

We can artificially introduce the dual map into the theory of Atkin-Lehner operators. For instance, one could write the Fricke involution $W_N$ as
\begin{displaymath}
	W_N f (z) = f \left( \begin{pmatrix}
		& - 1\\
		N & 
	\end{pmatrix}
	z \right) = f \left( \begin{pmatrix}
		& - 1\\
		N & 
	\end{pmatrix} w z^{- T} w \right) = f \left( \begin{pmatrix}
		1 & \\
		& N
	\end{pmatrix} z^{- T} \right).
\end{displaymath}
Though slightly cumbersome in rank 1, this approach leads to the right definition of Fricke involutions for $n > 2$.

Let $g \in \GL_n (\mathbb{R})$ such that
\begin{equation} \label{eq-fund}
	g^{- 1} \Gamma_0 (N) g = \Gamma_0 (N)^T. 
\end{equation}
Then the map $f (z) \mapsto f (gz^{- T})$ is an operator on the space of automorphic forms for $\Gamma_0 (N)$, which we call by definition an \emph{Atkin-Lehner operator}. 
As in the previous example, all Atkin-Lehner operators for $n = 2$ can be interpreted as above. 
More precisely, taking a matrix in the normaliser of $\Gamma_0 (N) \leq \SL_2(\R)$ and multiplying from the right by the non-trivial Weyl element gives a matrix $g$ satisfying \eqref{eq-fund}.

We now provide an example of Atkin-Lehner operators for all $n$.
The author was informed that Gergely Harcos has also, independently, found an example in the case $n = 3$.

\begin{definition}
	Let
	\begin{displaymath}
		A_N = N^{-1/n} \diag(1, \ldots, 1, N)
	\end{displaymath} 
	and define the \emph{Fricke involution} $W_N: L^2(X_n(N)) \longrightarrow L^2(X_n(N))$ as
	\begin{displaymath}
		W_N \phi(z) = \phi(A_N \cdot z^{-T}).
	\end{displaymath}
	We often also refer to the Fricke involution at the group level and denote
	\begin{displaymath}
		z' := A_N \cdot z^{-T}.
	\end{displaymath}
\end{definition}

It is easy to check that $A_N$ satisfies \eqref{eq-fund}. The operator defined above is obviously an involution and the expected properties hold.

\begin{lemma} \label{lemma-Fricke-invol}
	The Fricke involution $W_N$ preserves the space of cuspidal newforms and is self-adjoint. If $T_g$ is the Hecke operator associated to the coset $\Gamma_0(N) g \Gamma_0(N)$, where $(\det(g), N) = 1$, then
	\begin{displaymath}
		T_g W_N = W_N T_g^*.
	\end{displaymath}
	If an automorphic form $\phi$ has spectral parameters $(\mu_1, \ldots, \mu_n)$, then $W_N \phi$ has parameters $(-\mu_n, \ldots, -\mu_1)$.
\end{lemma}
\begin{proof}
	We first prove that $T_g W_N = W_N T_g^*$. 
	By a variant of the Smith normal form, we may assume that $g$ is diagonal and by a variant of the transposition anti-automorphism for $\Gamma_0 (N)$ (generalising Lemma 4.5.2 and Theorem 4.5.3 in \cite{miyake}), we may assume that there are matrices $\alpha_i$, $i = 1, \ldots, k$, for some $k$, such that
	\begin{displaymath}
		\Gamma_0 (N) g \Gamma_0 (N) = \bigcup_i \Gamma_0 (N) \alpha_i =
		\bigcup_i \alpha_i \Gamma_0 (N).
	\end{displaymath}
	Then by definition we have
	\begin{displaymath}
		T_g W_N f (z) = \sum_i W_N f (\alpha_i z) = \sum_i f (A_N \cdot
		\alpha_i^{- T} z^{- T}) = \sum_i f (\beta_i \cdot A_N \cdot z^{- T}) =
		W_N \sum_i f (\beta_i z),
	\end{displaymath}
	where $\beta_i = A_N \alpha_i^{- T} A_N^{- 1}$. 
	The proof is finished by showing that $\bigcup_i \Gamma_0 (N) \beta_i = \Gamma_0 (N) g^{- 1} \Gamma_0 (N)$, since this double coset corresponds to $T_g^{\ast}$ (s. \cite[Thm. 9.6.3]{goldfeld}). 
	Indeed,
	\begin{align*}
		\bigcup_i \Gamma_0 (N) \beta_i & = \bigcup_i \Gamma_0 (N) W_N
		\alpha_i^{- T} A_N^{- 1}\\
		& = \bigcup_i A_N \Gamma_0 (N)^T A_N^{- 1} A_N \alpha_i^{- T} A_N^{-
			1}\\
		& = A_N \left[ \bigcup_i \Gamma_0 (N) \alpha_i \right]^{- T} A_N^{-
			1}\\
		& = A_N \Gamma_0 (N)^T g^{- 1} \Gamma_0 (N)^T A_N^{- 1}\\
		& = \Gamma_0 (N) g^{- 1} \Gamma_0 (N) .
	\end{align*}
	Here we made use of fundamental property \eqref{eq-fund} of $A_N$ and of the fact that $g$ is diagonal, thus commuting with $A_N$.
	
	Next, we prove that $W_N$ is self-adjoint. 
	This can easily be seen by using a known fact about the dual forms for $\SL_n (\Z)$. 
	Namely, the map $f (z) \mapsto f (wz^{- T} w^{- 1})$, where $w$ is the long Weyl element, is self-adjoint (one can compute directly in explicit coordinates given in \cite{goldfeld}, Proposition 9.2.1 or Proposition 6.3.1). 
	We can interpret the Fricke involution as
	\begin{displaymath}
		W_N f (z) = f (m w z^{- T} w^{- 1}),
	\end{displaymath}
	where $m = A_N w^{- 1}$, that is, as the composition of the dualising map above with the left-action of $m$. 
	Since the measure on $\mathbb{H}^n$ is $\GL_n (\R)$-invariant, we can make the same explicit computations and change of coordinates as for the dualising map. Since $A_N$ is diagonal, we easily deduce the conclusion $W_N^{\ast} = W_N$. 
	Moreover, this interpretation of the Fricke involution and \cite[(45)]{BHMn} also prove the statement about the spectral parameters of $\phi$. 
	
	To prove cuspidality it is best to work adelically, though this can be reduced again to noting the relation between $W_N$ and the dualising map.
	Namely, the form $W_N \phi$ generates the contragredient of the representation generated by $\phi$, which is known to be cuspidal (see e.g. \cite[Prop. 3.3.4]{bump}). 
	From this perspective, it is also easy to see that $W_N \phi$ is a newform. 
	In the interest of brevity, we leave out the details of adelisation.\footnote{Based on our definition, the adelic formulation of these operators and their application in computing epsilon factors in functional equations has been worked out by Dalal and Gerbelli-Gauthier \cite[Sec.~3]{dalal-gerbelli}.}
\end{proof}

In this interpretation of Atkin-Lehner operators, the group structure coming from the normaliser is not obvious any more. 
Indeed, using \eqref{eq-fund}, we cannot even recover the identity for $n > 2$.
Finding an even more general definition proves difficult, since the available types of automorphisms on invertible matrices are scarce.

As explained in {\cite{mcd}}, all automorphisms in the case $n > 2$ are constructed out of inner automorphisms, radial automorphisms,
and the inverse-transpose automorphism. 
Inner automorphisms cannot contribute, since we have proved that the normaliser of $\Gamma_0 (N)$ is trivial; radial automorphisms are trivial in our context, since we consider only automorphic forms that are invariant under the centre of $\GL_n (\mathbb{R})$; and the
inverse-transpose automorphism is precisely the basis for the definition given in this note.

\subsubsection{Uniqueness of the Fricke involution} \label{sec-fricke-uniq}

The theory of Atkin-Lehner operators for $\Gamma_0 (N)$ shows some weaknesses already in the well-understood case $n = 2$. 
Indeed, one can only define Atkin-Lehner operators for divisors $M$ of the level $N$, such that $M$ and $N / M$ are coprime. 
More precisely, there are no operators induced by matrices with determinant equal to $M \mid N $, such that $(M, N / M) \neq 1$ (see \cite[p.\ 138]{AL}). 

This phenomenon creates difficulties in applications when considering powerful levels, as already noted in the historical context of the sup-norm problem. 
Here, we see that these difficulties only get more problematic in higher rank. 
In fact, the only Atkin-Lehner operator for $n > 2$, according to our definition, is the Fricke involution.

\begin{corollary} \label{prop-only-fricke}
	Let $g \in \GL_n^+ (\mathbb{Q})$ satisfy $g^{- 1} \Gamma_0 (N) g =
	\Gamma_0 (N)^T$. Then, after scaling by a suitable rational number, $g$ is
	integral, the last row and the last column of $g$ are divisible by $N$, and
	$\det (g) = N$. Equivalently, 
	\begin{displaymath}
		g \in \mathbb{Q}_{> 0} \cdot \Gamma_0 (N) \diag(1, \ldots, 1, N).
	\end{displaymath}
\end{corollary}

\begin{proof}
	This follows from Theorem \ref{thm-normaliser} on the triviality of the normaliser of $\Gamma_0(N)$.
	Indeed, let $g$ be as in the statement and recall the definition of the matrix $A_N$ that gives the Fricke involution.
	We have that
	\begin{displaymath}
		A_N \Gamma_0(N)^T A_N^{-1} = \Gamma_0(N),
	\end{displaymath}
	but also
	\begin{displaymath}
		A_N \Gamma_0(N)^T A_N^{-1} = A_N g^{-1} \Gamma_0(N) g A_N^{-1}.
	\end{displaymath}
	Thus, $g A_N^{-1}$ belongs to the normaliser of $\Gamma_0(N)$ and the result follows.
\end{proof}

\section{Reduction of the domain} \label{sec-red}

After studying generalised Atkin-Lehner operators, we showcase their main application in this section.
More precisely, we study fundamental domains for the action of these operators on $X_n(N)$.
Though very natural at a geometric level, we first note how this is relevant to the sup-norm problem.

For some $w \in \mathbb{H}^n$, the value of $\phi(w)$ is independent of which element in the orbit $\Gamma_0(N) \cdot w$ we choose instead of $w$. 
Similarly, the number and shape modulo $N$ of the matrices we are considering in the amplified pre-trace formula in Proposition \ref{prop-amp} is invariant under shifting by elements of $\Gamma_0(N)$, which would merely amount to conjugating $H(w, m, N)$.

Consider now the action of the Fricke involution $W_N(\phi)(w) = \phi(w')$. 
If $Y \subset \Gamma_0(N)$ is a subset, we denote by $Y'$ the image of $Y$ under the map $w \mapsto w'$. 
It is clear that we obtain a bound for a Hecke-Maaß form $\phi$ on $Y \cup Y'$ if we have a bound for both $\phi$ and $W_N(\phi)$ on the subset $Y$. 

Recall now that $W_N(\phi)$ has essentially the same properties as $\phi$ by Lemma \ref{lemma-Fricke-invol}.
Since the amplifiers, Proposition \ref{prop-amp} and Proposition \ref{prop-amp-uncond}, and the Fourier bound, Proposition \ref{prop-fourier}, apply similarly to both forms, we are free to choose any representative in
\begin{displaymath}
	\Gamma_0(N) w \cup \Gamma_0(N) w'
\end{displaymath}
when attacking the counting problem.\footnote{Indeed, the implied constant depending on $\mu$ in the amplifier is also of the same size, as the computation of spectral parameters in Lemma \ref{lemma-Fricke-invol} shows.} 

In this section we propose a system for making this selection of representative.
In other words, we construct an approximate fundamental domain for the action of $\Gamma_0(N)$ and the Fricke involution, at least on $\Omega_N^\natural$. 
It can be seen as a reduction theory with level structure, for which we often use the shorter term Fricke reduction.

\subsection{A reinterpretation of the Harcos--Templier reduction}\label{sec-HT-red}

To motivate our approach, we first revisit the case of $\PGL(2)$, assuming that $N$ is prime.  
Here, the elegant and simple idea of Harcos--Templier \cite{HT2} is to choose a representative
\begin{displaymath}
	z \in \Gamma_0(N) w \cup \Gamma_0(n) w'
\end{displaymath}
such that $\Im z$ is maximal when viewing $z$ in the upper-half plane.
Using the Iwasawa coordinates $z = n(x) a(y)$, where $a(y) = d \cdot \diag(y,1)$, this amounts to minimising $\norm{e_2}_z$.
What this achieves is that $y \gg 1/N$ and, in the upper-half plane, $|cz + d| \geq 1/\sqrt{N}$ for all $(c,d) \in \Z^2$ (see Lemma 1 in \cite{HT2}).
After normalisation, the latter condition is equivalent to saying the shortest length in $L_z$ is bounded from below by $d / \sqrt{N}$.
Note that $d = 1/\sqrt{y}$.

This last bound shows that $y$ controls the balancedness of $L_z$, which plays a crucial role in the counting problem.
For $y$ close to $1/N$, the lattice $L_z$ is balanced and this is tractable, yet for larger $y$ they resort to applying a Fourier bound.

We now look closer at this reduction process and compare the two lattices $L_{w}$ and $L_{w'}$.
Define
\begin{displaymath}
	\alpha(w) = \min \{\norm{e_2}_{\gamma z} \mid \gamma \in \Gamma_0(N) \} 
\end{displaymath}
and similarly for $\alpha(w')$.
Now $e_2 \cdot \Gamma_0(N)$ does not exhaust all primitive vectors in $\Z^2$ and thus $\alpha(w)$ does not necessarily give the first successive minimum of $L_w$.
However, since $N$ is prime, we only lack primitive vectors of the form $(a, b)$, where $\gcd(a, N) = 1$.
In this case $(Nb, a) \in \Z^2$ is another primitive vector.
We may thus find $\gamma \in \Gamma_0(N)$ with last row equal to $(Nb, a)$.

Consider, for simplicity, the old definition of the Fricke involution
\begin{displaymath}
	w' = N^{-1/2} \begin{pmatrix}
		0 & -1 \\
		N & 0
	\end{pmatrix} w.
\end{displaymath}
A quick computation shows that
\begin{displaymath}
	\norm{e_2}_{\gamma w'} = \norm{(Nb, a)}_{w'} = \sqrt{N} \norm{(a, b)}_w.
\end{displaymath}
Putting everything together, the shortest length in $L_w$ is either $\alpha(w)$ or $\alpha(w')/\sqrt{N}$.
By symmetry, the shortest length in $L_{w'}$ is either $\alpha(w')$ or $\alpha(w)/\sqrt{N}$.

To choose a representative $z$, what Harcos--Templier do is to potentially apply the Fricke involution to achieve that $\alpha(z) \leq \alpha(z')$ and then to shift $z$ by $\Gamma_0(N)$ to achieve $\alpha(z) = d$, where $d = 1/\sqrt{\Im z}$.
We now look at what information this choice implies about the lattices above, and we distinguish a few possibilities.

If $\alpha(z) = \norm{e_2}_z = d$ is the shortest length in $L_z$, then we have achieved the first step in Gauß--Lagrange reduction, and we may continue to obtain $y \geq \sqrt{3}/2$. 
For this conclusion we do not use that $\alpha(z) \leq \alpha(z')$.
In fact, since this lower bound on $y$ is so strong, it seems useful to prioritise using the Fricke involution to arrive at this situation.
Thus, in the case that $\alpha(z')$ is the shortest length in $L_{z'}$, we could apply the Fricke involution and obtain the same conclusion.

We are left with the case where the shortest length in $L_z$ is $\alpha(z')/\sqrt{N}$ and the shortest in $L_{z'}$ is $\alpha(z)/\sqrt{N}$.
We assume the Harcos--Templier conditions and investigate what happens if we assume that our original point reduces to a compactum $\Omega$. 
Interpreting this in terms of balancedness of lattices but not knowing whether the Fricke involution was applied or not, we must distinguish two cases.

First, if we applied the involution and $z'$ reduces to $\Omega$, meaning that $L_{z'}$ is balanced, then its shortest length $\alpha(z)/\sqrt{N}$ is approximately $1$, implying that $d = 1/\sqrt{y} \asymp_\Omega \sqrt{N}$.
Thus, we have the case $y = 1/N$, where the counting results work best.
As noted before, this also implies that $L_z$ is balanced.
Indeed, the shortest length in $L_z$ is $\alpha(z')/\sqrt{N} \geq \alpha(z)/\sqrt{N} \asymp 1$, by assumption.

Otherwise, if $L_z$ is the balanced lattice and we cannot assume anything about $L_{z'}$, notice that the crucial parameter $y$ still controls the balancedness of $L_{z'}$.
Thus, either $L_{z'}$ is rather balanced, or $y$ is much larger than $1/N$.
This dichotomy can be successfully used in solving the sup-norm problem, as long as the Fourier bound is effective for $y \gg N^{-1 + \delta}$ for small enough $\delta$.

The case of $y$ much larger than $1/N$ seems possible, since no contradictions can be found with our assumptions.
In fact, $L_z$ being balanced implies the trivial bound $\alpha(z) \gg 1$, translating only to $y \ll 1$. 
For instance, this shows that the shortest length in $L_{z'}$ is bounded from below by roughly $1/\sqrt{N}$.

The considerations above show that the preimage of a compact set $\Omega$ under the natural projection map $X_2(N) \to X_2(1)$ can be covered by $Y \cup Y'$, where
\begin{displaymath}
	Y \subset \{ z = x + i y \mid x \in [-1/2, 1/2], y \gg 1 / N \}.
\end{displaymath}
Harcos--Templier show the same for the full preimage of $X_2(1)$, yet the inclusion above seems sharp even for a compact set.

\subsection{Two lattices} \label{sec-two-lattices}

We now extend the ideas above to higher rank, which is significantly more complicated. 
Throughout the following sections we assume that $N$ is a prime. 

Recall that for $z \in \SL_n(\R)$ we write 
\begin{displaymath}
z' := A_N w^{-T} = N^{-1/n} \diag(1, \ldots, 1, N) z^{-T}.
\end{displaymath} 
We consider the lattices $L_z$ and $L_{z'}$ in the notation and terminology established in Section \ref{sec-prelims}. 
Note that both lattices have determinant $1$. 
We define the sets
\begin{align*}
	A(z) &= \{ \norm{e_n}_{\gamma z} \mid \gamma \in \Gamma_0(N) \}, \\
	B(z) &= \{ \norm{e_2 \wedge \cdots \wedge e_n}_{\gamma z} \mid \gamma \in \Gamma_0(N) \}.
\end{align*}
In the following paragraphs we show how the union of $A(z), B(z), A(z'), B(z')$ provides the lengths of all primitive vectors in $L_z$, $L_{z'}$, and their duals.

First, we claim that the union of lengths
\begin{displaymath}
	\{ \norm{e_n}_{\gamma z} \mid \gamma \in \Gamma_0(N) \} \cup \{ \norm{e_2 \wedge \cdots \wedge e_n}_{\gamma z'} \mid \gamma \in \Gamma_0(N) \}
\end{displaymath}
exhausts the lengths of all primitive vectors in $L_z$. For this we use the fact that any primitive vector in $\Z^n$ is the last row (in fact, any row or any column) of some matrix in $\SL_n(\Z)$. Consequently, the vectors $e_n \gamma$ give all primitive vectors in $N \Z \times \cdots \times N \Z \times \Z$ in the lattice $L_z$. 

For the second set, note using Lemma \ref{lemma-isometric-lats} that
\begin{equation} \label{eq-norms-z-1}
	\norm{e_2 \wedge \cdots \wedge e_n}_{\gamma z'} = \norm{e_1}_{\gamma^{-T} A_N^{-1} z} = N^{1/n} \norm{(a_1, \ldots, a_n)}_z,
\end{equation}
where $(a_1, \ldots, a_{n-1}, Na_n)$ is the top row of $\gamma^{-T}$. We prove in Lemma \ref{lemma-matrix-completion} below that we obtain this way all primitive vectors $(a_1,\ldots, a_n)$ in $L_z$, for which 
\begin{displaymath}
	\gcd(a_1, \ldots, a_{n-1}, N) = 1.
\end{displaymath}
Since $N$ is prime, the greatest common divisor of $a_1, \ldots, a_{n-1}$ and $N$ can only be $1$ or $N$.
Thus, considering the paragraph above, we have exhausted all primitive vectors in~$L_z$.

\begin{lemma} \label{lemma-matrix-completion}
	For $N$ prime, if $v = (a_1, \ldots, a_{n-1}, Na_n) \in \Z^n$ is a primitive vector, then there is $\gamma \in \Gamma_0(N)$ such that $v$ is the first row of $\gamma^T$.
\end{lemma}
\begin{proof}
	Let $g \in \SL_n(\Z)$ be any matrix with first row $v$. Multiplying $g$ from the left by block matrices of the form
	\begin{displaymath}
		\begin{pmatrix}
			1 & \\
			& h
		\end{pmatrix},
	\end{displaymath}
	where $h \in \SL_{n-1}(\Z)$, leaves the first row invariant. We shall inductively apply such row operations on $g$ to make its last column be of the form $(c_1, \ldots, c_n)$, where $c_1 = N a_n$ and $N$ divides $c_1, \ldots, c_{n-1}$. 
	
	Indeed, if $N \mid c_i$ for any $i \in \{2, \ldots, n\}$, then we can permute rows to assume that $N \mid c_2$. 
	Otherwise, we can assume that $\gcd(N, c_3) = 1$. Let $\bar{c}_3$ be any representative of the inverse of $c_3$ modulo $N$. 
	Bézout's lemma provides a matrix $h' \in \SL_2(\Z)$ with top row $(N, \bar{c}_3)$. 
	Using $h$ of the form
	\begin{displaymath}
		h = \begin{pmatrix}
			h' & \\
			& 1_{n-3}
		\end{pmatrix},
	\end{displaymath}
	as above, we may now assume that $c_2 \equiv 1$ modulo $N$. Another transformation of the same type, where $h'$ now has top row $(1, -\bar{c}_3)$, allows us to assume that $N \mid c_2$. We conclude by induction.
\end{proof}

Next, the union of lengths
\begin{displaymath}
	\{ \norm{e_2 \wedge \cdots \wedge e_n}_{\gamma z} \mid \gamma \in \Gamma_0(N) \} \cup \{ \norm{e_n}_{\gamma z'} \mid \gamma \in \Gamma_0(N) \}
\end{displaymath}
exhausts the lengths of all primitive vectors in $L^\ast_{z} = L_{z^{-T}}$. Indeed, Lemma \ref{lemma-isometric-lats} gives that 
\begin{equation}  \label{eq-norms-ztilde-1}
	\norm{e_2 \wedge \cdots \wedge e_n}_{\gamma z} = \norm{e_1}_{\gamma^{-T} z^{-T}} = \norm{(a_1,\ldots a_{n-1},N a_n)}_{z^{-T}},
\end{equation}
where $(a_1,\ldots,Na_n)$ is the first row of $\gamma^{-T}$. As above, we obtain this way all primitive vectors in $\Z^{n-1} \times N \Z$ in the lattice $L_{z^{-T}}$. Furthermore, 
\begin{equation} \label{eq-norms-ztilde-2}
	\norm{e_n}_{\gamma A_N z^{-T}} = N^{1-1/n} \norm{(a_1,\ldots,a_n)}_{z^{-T}},
\end{equation}
for $(a_1,\ldots,a_n)$ primitive with $\gcd(a_n, N)=1$. Since $N$ is prime, this shows the claim.

The above considerations are collected for an overview in Table \ref{table}.
\begin{table}[h]
	\renewcommand{\arraystretch}{1.5}
	\centering
	\begin{tabular}{c||c|c}
	$L_z$ & $A(z)$ & $N^{-1/n} \cdot B(z')$   \\
	\hline
	$L_z^\ast$ & $N^{-1+1/n} \cdot A(z')$ & $B(z)$   \\
	\hline
	$L_{z'}$ & $A(z')$ & $N^{-1/n} \cdot B(z)$   \\
	\hline
	$L_{z'}^\ast$ & $N^{-1+1/n} \cdot A(z)$ & $B(z')$
	\end{tabular}
	\caption{Lattices and sets of lengths of primitive vectors.}
	\label{table}
\end{table}
Each row corresponds to a lattice and the union of the two sets in that row is the set of the lengths of all primitive vectors in the corresponding lattice. By multiplication of a set by a scalar we mean multiplication of each element in the set by the given scalar.
We use here that $z \mapsto z'$ is an involution on unimodular lattices. 

\subsection{Fricke reduction} \label{sec-fricke-red}

Let us consider minima of the lattices in the previous section. Write
\begin{displaymath}
	\alpha(z) = \min A(z), \quad \beta(z) = \min B(z).
\end{displaymath}
As in Table \ref{table}, the minimal non-zero length in the lattice $L_z$ is found either in $A(z)$, equal in this case to $\alpha(z)$, or in $B(z')$, equal to $N^{-1/n} \beta(z')$.

More generally, let $x$ be any of the letters $\alpha$ or $\beta$. Let $L$ be any of the lattices $L_z$, $L_{z'}$, $L_z^\ast$, $L_{z'}^\ast$. Then the minimal length in $L$ is an $x$-expression if it is of the form $N^\eta x(w)$, where $\eta$ is a non-positive number and $w$ is either $z$ or $z'$. From Table \ref{table} and the discussion of that section, we see that there are only two possibilities for each lattice, namely a unique $\alpha$-expression or a unique $\beta$-expression.

\begin{definition}
	Let $X$ and $Y$ denote the \emph{Greek} letters $A$ or $B$, and analogously for their lowercase variants. We say that $z \in \mathcal{L}(X, Y)$ if the smallest length in $L_z$ is the unique $x$-expression and the smallest length in $L_{z}^\ast$ is the unique $y$-expression. Similarly, $z \in \mathcal{L}'(X, Y)$ if the smallest lengths in $L_{z'}$ and $L_{z'}^\ast$ are the $x$-expression and the $y$-expression, respectively.
\end{definition}

\begin{example}
	If $z \in \mathcal{L}(B, A)$, then the smallest length in $L_z$ is given by $N^{-1/n} \beta(z')$ and the smallest length in $L_{z}^\ast$ is given by $N^{-1+1/n} \alpha(z')$.
\end{example}

Let $w \in \mathbb{H}$. For the study of the sup-norm and our counting problem, we are allowed to choose any conjugate of $z$ in the orbit $\Gamma_0(N) \cdot w$ and also switch between $z$ and $z'$, as explained at the beginning of Section \ref{sec-red}. 
Now, it is clear by construction that every $w$ is contained in some $\mathcal{L}(X,Y)$. We then make the choice of conjugate to obtain a well-positioned $z$, where we have control over its successive minima and Iwasawa coordinates, based on which set $\mathcal{L}(X,Y)$ contains $z$. 
The guiding principle is that $y$-coordinates of $z$ should be bounded from below and as large as possible, similar to classical reduction theory.
If this is not possible, we try to obtain some inbalancedness information on the relevant lattices.

We first begin with some preliminary general observations and then construct a domain in the following proposition, under a balancedness assumption.

\subsubsection{Case I} \label{sec-case-I}
Let 
\begin{displaymath}
	z \in \bigcup_{X \in \{A, B\}} \mathcal{L}(A, X) \cup \mathcal{L}'(A, X).
\end{displaymath}
By switching between $z$ and $z'$ if needed, we can assume that $z \in \mathcal{L}(A,X)$, for some $X \in \{A, B\}$. In this case, the minimal length in $L_z$ is $\alpha(z)$. Shifting $z$ by $\gamma \in \Gamma_0(N)$ if needed, we assume that $\alpha(z) = \norm{e_n}_z$. In Iwasawa coordinates $z = n(x) \cdot a(y)$ as in Section \ref{sec-Iwasawa}, we have $\alpha(z) = d$.

Let $\gamma$ be of the form
\begin{displaymath}
	\gamma = \begin{pmatrix}
		h & \\
		& 1
	\end{pmatrix} \in \Gamma_0(N),
\end{displaymath}
where $h \in \SL_{n-1}(\Z)$. Note that $e_n \cdot \gamma = e_n$, so we can make the same assumptions about $\gamma z$ as about $z$ above. As in Remark \ref{rem-reduction-n-1}, shifting by $\gamma$ as above if needed, we may now additionally assume that $z = n(x) a(y)$ satisfies $y_i \geq \sqrt{3}/2$ for $i = 2, \ldots, n-1$.

By Lemma \ref{lemma-Mahler}, if $\lambda_1$ and $\lambda_2$ are the first two successive minima of $L_z$, then the shortest length $l$ in $\bigwedge^2 L_z$ satisfies
\begin{displaymath}
	l \asymp_{n} \lambda_1 \cdot \lambda_2.
\end{displaymath}
In particular, $l \gg \lambda_1^2$. This implies that (recall the calculation in Remark \ref{rem-iwasawa-calculations})
\begin{displaymath}
	d^2 y_1 = \norm{e_{n-1} \wedge e_n}_z \gg \alpha(z)^2 = d^2.
\end{displaymath}
We deduce that $y_1 \gg_n 1$.\footnote{This can be viewed as a soft version of Hermite reduction, that is, reduction to a Siegel set. Indeed, here we also take the last row to be the shortest vector and then use induction, as in the classical proof of reduction.}

\subsubsection{Case II}
Let
\begin{displaymath}
	z \in [\mathcal{L}(B,A) \cap \mathcal{L}'(B,A)] \cup [\mathcal{L}(B,B) \cap \mathcal{L}'(B,A)].
\end{displaymath} 
Shifting $z = n(x) a(y)$ by a suitable $\gamma \in \Gamma_0(N)$ as in Case I, we assume that $\alpha(z) = \norm{e_n}_z = d$ and that $y_i \gg 1$ for $i = 2, \ldots, {n-1}$. 

Note now that the minimal length in $L_{z'}^\ast$ is $N^{-1+1/n} \alpha(z)$. 
Note also that $L_{z'}^\ast$ is the lattice corresponding to
\begin{displaymath}
	z'^{-T} = A_N^{-1} z.
\end{displaymath}
We now compute that
\begin{displaymath}
	\norm{e_{n-1} \wedge e_n}_{z'^{-T}} = \norm{e_{n-1} \wedge e_n}_{A_N^{-1} z} = N^{-1 + 2/n} d^2 y_1.
\end{displaymath}
Using Lemma \ref{lemma-Mahler}, we deduce that the minimal length $l$ in $\bigwedge^2 L_{z'}^\ast$ satisfies $l \gg \mu_1^2$, where $\mu_1$ is the first successive minimum of $L_{z'}^\ast$. 
Putting everything together we arrive at
\begin{displaymath}
	N^{-1 + 2/n} d^2 y_1 \gg N^{-2 + 2/n} d^2,
\end{displaymath}
which implies that $y_1 \gg_n N^{-1}$.
Notice that the Fricke involution has not been used.

\subsubsection{Case III}
Let
\begin{displaymath}
	z \in [\mathcal{L}(B,A) \cap \mathcal{L}'(B,B)] \cup [\mathcal{L}(B,B) \cap \mathcal{L}'(B,B)].
\end{displaymath}
Then the minimal length in $L^\ast_{z'}$ is given by $b(z')$ and the minimal length in $L_z$ is $N^{-1/n} b(z')$.
By Minkowski's theorem, more precisely equation \eqref{eq-lambda-1}, applied to $L_{z'}^\ast$, we find that $b(z') \ll_n 1$. 
This now implies that the minimal length in $L_z$ is $N^{-1/n} b(z') \ll N^{-1/n}$.
Notice that we again do not need the Fricke involution for this conclusion.

\subsubsection{Fricke reduction of points that reduce to a compactum}
We summarise the cases described above in the context of points $w$ that reduce to a fixed compact set $\Omega \subset \mathbb{H}$ and we produce a well-positioned conjugate $z$.

\begin{proposition} \label{prop-fricke-red}
	Let $w \in \mathbb{H}^n$, let $\Omega \subset \mathbb{H}^n$ be a compact set and assume that $w$ reduces to~$\Omega$.
	For $N \gg_\Omega 1$ prime, large enough, only Case I and Case II are valid.
	There is
	\begin{displaymath}
		z \in \{ \gamma w \mid \gamma \in \Gamma_0(N) \} \cup \{ \gamma w' \mid \gamma \in \Gamma_0(N) \},
	\end{displaymath}
	where $w' = A_N w^{-T}$, with Iwasawa coordinates $z = n(x) a(y)$,
	\begin{displaymath}
		a(y) = d \cdot \begin{pmatrix}
			y_1 \cdots y_{n-1} & &  & \\
			& \ddots &  & \\
			& & y_1 & \\
			& & & 1
		\end{pmatrix},
	\end{displaymath} 
	which satisfies:
	\begin{enumerate}
		\item if $w$ is in Case I, then
		\begin{displaymath}
			y_i \gg 1 \text{ for all } i = 1, \ldots, n-1;
		\end{displaymath}

		\item if $w$ is in Case II, then we can arrange that $d = \alpha(z)$, we have the lower bounds
		\begin{displaymath}
			y_1 \gg \frac{1}{N} \quad \text{and} \quad y_i \gg 1 \text{ for }i = 2, \ldots, n-1,
		\end{displaymath} 
		and either
		\begin{enumerate}
			\item $w$ belongs to a set of volume $\ll N^{(n-1)(1-1/2n^4)}$ and we have 
			\begin{displaymath}
				y_1 \gg N^{-1 + 1/2n^4} \text{ or } y_i \gg N^{1/2n^4} \text{ for some } 2 \leq i \leq n-1;
			\end{displaymath}
			\item or $z, w \in \mathcal{L}(B,A) \cap \mathcal{L}'(B,A)$, the lattice $L_z$ is $\Omega$-balanced, we have 
			\begin{displaymath}
				\lambda_1(L_{z'}) \gg N^{-1/2n}
			\end{displaymath}
			and we have the upper bounds
			\begin{displaymath}
				dy_1 \cdots y_k \ll N^{-1/2n} \text{ for any } 1 \leq k \leq n-1.
			\end{displaymath}
		\end{enumerate}
	\end{enumerate}
	All implied constants may depend on $\Omega$ and $n$.
\end{proposition}
\begin{proof}
	In \textbf{Case I}, we have already shown how to find $z$ as in the statement such that $y_i \gg_n 1$ for all $i = 1, \ldots, n-1$.
	
	\textbf{Case II} is more delicate and we first consider the possibility that
	\begin{displaymath}
		w \in \mathcal{L}(B,B) \cap \mathcal{L}'(B,A).
	\end{displaymath} 
	Here, as explained above, we may find $z=\gamma w$ with $\gamma \in \Gamma_0(N)$ with Iwasawa coordinates satisfying that $d = \alpha(z)$, $y_1 \gg N^{-1}$, and $y_i \gg 1$ for $i \geq 2$.
	Observe now that the shortest vector in $L_z^\ast$ is given by $\beta(z)$.
	Since we did not apply the Fricke involution, $z$ reduces to $\Omega$.
	By Lemma \ref{lemma-balancedness}, $L_z^\ast$ is a balanced lattice and therefore $\beta(z) \asymp_\Omega 1$.

	On the other hand, the shortest vector in $L_{z'}$ is given by $N^{-1/n} \beta(z)$ and this is now of size $N^{-1/n}$.
	By Minkowski's theorem, we have $\lambda_1(L_{z'}) \lambda_n(L_{z'})^{n-1} \gg 1$ and therefore $\lambda_n(L_{z'}) \gg N^{1/n(n-1)}$.
	By a standard transference theorem, putting together Lemma \ref{lemma-isometric-lats} and Lemma \ref{lemma-Mahler}, we find that
	\begin{displaymath}
		\lambda_1(L_{z'}^\ast) \ll \lambda_n(L_{z'})^{-1} \ll N^{-1/n(n-1)}.
	\end{displaymath}
	Since the shortest vector in $L_{z'}^\ast$ is given by $N^{-1+1/n} \alpha(z)$, our assumptions imply that $d \ll N^{1-1/n-\eta_n}$, where $\eta_n = 1/n(n-1)$.

	Since $d^n y_1^{n-1} \cdots y_{n-1} = 1$ by assumption, we deduce from the upper bound on $d$ that $y_1 \gg N^{-1 + \eta_n'}$ or $y_i \gg N^{\eta_n'}$ for some $i \geq 2$, where we may take $\eta_n' = \eta_n/(n-1)^2$.
	We may bound now the volume of the Siegel-like set with such lower bounds on the $y$-coordinates (and $x$-coordinates uniformly bounded, as usual).
	Recall (see \cite[Prop.~1.5.3]{goldfeld}) that the volume form for $\mathbb{H}_n$ in the $y$-coordinates is
	\begin{displaymath}
		\prod_{k=1}^{n-1} y_k^{-k(n-k)-1} \, dy_k.
	\end{displaymath}
	The largest exponent appearing is $-n$ and therefore the bounds above show that the volume is bounded by $N^{(n-1)(1-\eta_n')}$.

	We now consider the crucial case of
	\begin{displaymath}
		w \in \mathcal{L}(B,A) \cap \mathcal{L}'(B,A).
	\end{displaymath}
	This case is invariant under the Fricke involution, and we apply it if necessary to obtain $z$ such that $\alpha(z) \leq \alpha(z')$.
	Shifting by $\Gamma_0(N)$, we may also assume that the $y$-coordinates of $z$ satisfy $d = \alpha(z)$, $y_1 \gg_n N^{-1}$, and $y_i \gg 1$ for $i = 2, \ldots, n-1$, as we generally do with Case II.

	Since we potentially applied the Fricke involution, we must now consider two possibilities.
	First, if we did apply the involution and it is $z'$ that reduces to $\Omega$, then $L_{z'}$ and $L_{z'}^\ast$ are $\Omega$-balanced.
	Since $z \in \mathcal{L}'(B,A)$, the minimal length in $L^\ast_{z'}$ is $N^{-1 + 1/n} \alpha(z)$ and we deduce that $\alpha(z) = d \asymp N^{1-1/n}$.
	We see that
	\begin{displaymath}
		d^{-n} = y_1^{n-1} \cdots y_{n-1} \ll N^{-(n-1)}.
	\end{displaymath}
	Combining this with the bounds above for the $y$-coordinates, we deduce that $y_1 \asymp N^{-1}$ and $y_i \asymp 1$ for $i = 2, \ldots, n-1$, where the implicit constants depend on $\Omega$.
	Additionally, $z \in \mathcal{L}(B,A)$ implies that the shortest length in $L_z^\ast$ is 
	\begin{displaymath}
		N^{-1 + 1/n} \alpha(z') \geq N^{1+1/n} \alpha(z) \gg 1.
	\end{displaymath}
	As such, $L_z^\ast$ is a balanced lattice and, by Lemma \ref{lemma-balancedness}, the same is true about $L_z$.

	Otherwise, we have that $z$ reduces to $\Omega$, implying that $L_z$ is balanced.
	In this case, we cannot deduce that $L_{z'}$ is balanced, but any imbalancedness translates into lower bounds for $y$-coordinates, as above.

	Indeed, by assumption, the shortest vector in $L_{z'}$ has length $N^{-1/n} \beta(z)$.
	If $N^{-1/n} \beta(z) \leq N^{-1/2n}$, then the same arguments as above show that 
	\begin{displaymath}
		d \ll N^{1- 1/n - 1/2n(n-1)}, \quad y_1 \gg N^{-1 + 1/2n(n-1)^3}, \quad y_i \gg N^{1/2n(n-1)^3} \text{ for all } i \geq 2.
	\end{displaymath}
	The volume of such $z$ is bounded by $N^{(n-1)(1-1/2n^4)}$.

	Otherwise, assume that $N^{-1/n} \beta(z) \geq N^{-1/2n}$. 
	By minimality of $\beta(z)$, we have the bound
	\begin{displaymath}
		\beta(z) \leq \norm{e_2 \wedge \cdots \wedge e_n}_z = (d y_1 \cdots y_{n-1})^{-1},
	\end{displaymath}
	using the computation of Remark \ref{rem-iwasawa-calculations}.
	Therefore,
	\begin{displaymath}
		d y_1 \cdots y_{n-1} \leq N^{1/2n - 1/n} = N^{-1/2n}.
	\end{displaymath}
	Since $y_i \gg 1$ for all $i \geq 2$, this also implies that
	\begin{displaymath}
		d y_1 \cdots y_k \ll N^{- 1/2n}
	\end{displaymath}
	for all $1 \leq k \leq n-1$.

	Finally, by Lemma \ref{lemma-red-to-cpt}, we eliminate \textbf{Case III}, where we have shown that the minimal length in $L_w$ is $\ll N^{-1/n}$.
	Indeed, since Case III does not involve the Fricke involution, we preserve the assumption that $L_w$ is $\Omega$-balanced and arrive at a contradiction.
\end{proof}

\begin{remark}
	What we call the bulk $\Omega_N$ in the introduction of this paper is given here by case \textbf{(2.b)}.
	The cuspidal regions that we remove are those corresponding to \textbf{(1)} and \textbf{(2.a)}, as noted in \eqref{eq:exceptional}.
\end{remark}

\section{Counting matrices}

\subsection{An overview}

When applying the amplified pre-trace formula, e.g.\ Proposition~\ref{prop-amp-uncond}, we arrive at the problem of counting matrices in $H(z, m, N)$. 
We give a brief overview of the counting strategy in the simplest case of $n = 2$.
The perspective taken in this paper is new even in this case.
We recall some ideas already introduced in Section~\ref{sec-intro-lats}.

Let $z \in \SL_2(\R)$, for which we assume the Iwasawa form
\begin{displaymath}
	z = \begin{pmatrix}
		\sqrt y & x / \sqrt{y} \\
		0 & 1/\sqrt{y}
	\end{pmatrix}
	= n(x) \cdot d \cdot \diag(y, 1),
\end{displaymath}
where $d = 1/\sqrt{y}$, and let $\gamma \in H(z, m , N)$. 
The bound
\begin{equation} \label{eq-counting-prob-2}
	z^{-1} \gamma z = O(m^{1/2})
\end{equation}
implies the conditions
\begin{displaymath}
	e_i \cdot \gamma z \in B(m^{1/2} \norm{e_i \cdot z})
\end{displaymath}
for $i = 1, 2$, where $B(r)$ is a Euclidean ball of radius $O(r)$ around $0$.

We assume now that $z$ lies in bulk $\Omega_N$, i.e. it satisfies the conditions in case \textbf{(2.b)} of Proposition \ref{prop-fricke-red}.
In this context, we have that $L_z$ is a balanced lattice, and $1/N \ll y \ll 1/\sqrt{N}$.
Importantly, the shortest length in $L_{z'}$ is given by $d/\sqrt{N}$ (note either that rank 2 lattices are isometric to their duals or recall the discussion in Section \ref{sec-HT-red}).
Let 
\begin{displaymath}
	z_N = \diag(N, 1) \cdot z,
\end{displaymath}
which defines a sublattice of index $N$ of $L_z$.
In fact, using the classical definition of the Fricke involution and observing that
\begin{displaymath}
	\begin{pmatrix}
		0 & -1 \\
		N & 0
	\end{pmatrix} = \begin{pmatrix}
		0 & -1 \\
		1 & 0
	\end{pmatrix}
	\cdot 
	\begin{pmatrix}
		N & 0 \\
		0 & 1
	\end{pmatrix},
\end{displaymath}
we easily see that $L_{z_N}$ is equal to $N^{1/2} L_{z'}$.
As such, the shortest length in $L_{z_N}$ is given by $d = \norm{e_2}_z = \norm{e_2}_{z_N}$.

This is helpful since we now count the possibilities for $e_2 \cdot \gamma$, a vector in the sublattice $N\Z \times \Z$.
We do this by applying Lemma \ref{lemma-lattice-count}, which counts lattice points in balls.
Since the radius of the ball is given by $\sqrt{m} \cdot \norm{e_2}_z$ and $\norm{e_2}_z = d$ is the smallest length in $L_{z_n}$, the bound we obtain is 
\begin{displaymath}
	1 + \frac{d \sqrt{m}}{d} + \frac{d^2 m}{d^2} \ll m,
\end{displaymath} 
by bounding $\lambda_1 \lambda_2 \gg \lambda_1^2$.

For $e_1 \cdot \gamma$, we notice that $\norm{e_1}_z$ is equal to $y + x/y \asymp 1/N + N x$.
Unfortunately, if $z$ is a balanced lattice, one can compute that we must have a bound $x \gg 1/\sqrt{N}$.
Thus, the norm above can be rather large.
Even though $L_z$ is balanced, the size of the ball would give a hopelessly large bound.

Fortunately, we notice that
\begin{displaymath}
	e_1 \cdot z - x e_2 \cdot z = (\sqrt y, x / \sqrt y) - x (0, 1/\sqrt y) = (\sqrt y, 0),
\end{displaymath}
by the Iwasawa decomposition or the Gram-Schmidt process.
The conditions above can be combined to show that
\begin{displaymath}
	e_1 \cdot \gamma z - x e_2 \cdot \gamma z \in B(m^{1/2} \norm{(\sqrt y, 0)}).
\end{displaymath}
Since $y \ll 1/\sqrt{N}$, we see that if $m \ll N^{1/2-\varepsilon}$, the ball we obtain has a small radius of size $o(1)$.
Since $L_z$ is a balanced lattice, we can only have at most one lattice point in such a small ball, regardless of its centre.
For every vector $e_2 \cdot \gamma$ fixed as above, this leaves at most one possibility for $e_1 \cdot \gamma$.
Therefore, the second row $e_2 \cdot \gamma$ already fixes the whole matrix $\gamma$.

This strategy gives a bound
\begin{displaymath}
	\# \bigcup_{l=1}^{m} H(z, l, N) \ll m
\end{displaymath}
if $m$ is small enough in terms of $N$.
A glance at Proposition \ref{prop-amp-uncond} shows that this bound is insufficient to obtain a saving when averaging over square determinants $l = p^2 q^2$ and thus $m = L^4$, in the notation of the proposition.

To refine the process above, we only partially fix the second row of $\gamma$.
This seems difficult to do in standard coordinates, that is, working with the exact entries of $\gamma$.
Instead, we choose a reduced basis, $v_1$ and $v_2$, for the lattice $L_{z_N}$.
An upshot of Fricke reduction is that we can choose $v_2 = e_2 \cdot z$ (we already noticed above that $e_2 \cdot z$ is a shortest vector in $L_{z_N}$).

We now write $e_2 \cdot \gamma z \in L_{z_N}$ in coordinates using $v_1$ and $v_2$.
By our conditions and the balancedness of the lattice, the coefficients for both basis vectors are bounded by $\sqrt{m}$.
In a first step, we only choose the coefficient of $v_1$, giving us $\sqrt{m}$ possibilities.

We now ask how many matrices $\gamma$ have such a coefficient.
For two such matrices $\gamma_1, \gamma_2$, the difference $\gamma_1 - \gamma_2$ would have last row equal to $c \cdot e_2$ with $c \ll \sqrt m$.
It would also satisfy \eqref{eq-counting-prob-2}.
These two observations imply that the strategy above applies to this difference.
The principle that the last row fixes the matrix now gives that $\gamma_1 - \gamma_2 = c \cdot \id_2$.

Applying the determinant to $\gamma_1 = \gamma_2 + c \cdot \id_2$ and assuming that $\gamma_1$ has a square determinant imply that $-c$ gives a solution to
\begin{displaymath}
	\chi_{\gamma_2}(X) = Y^2.
\end{displaymath}
We employ a theorem of Heath-Brown to count solutions to such equations and obtain adequate bounds for the amplified pre-trace formula in the non-degenerate case.

The degenerate case is precisely when the characteristic polynomial of $\gamma_2$ is a square.
This means that $\gamma_2$ is a parabolic matrix and therefore fixes a cusp.
For $\Gamma_0(N)$ with $N$ prime, there are two such cusps and these are conjugated by the Fricke involution.
This allows us to assume that $\gamma_2$ fixes the cusp at infinity and is therefore an upper triangular matrix, up to conjugation.
The strategy above can be adapted slightly for us to apply, again,  the principle that the last row determines the matrix.
In this case, the last row is the same as that of a multiple of the identity matrix and we are done.

\subsection{The iterative strategy} \label{sec-counting-strategy}

In this section we generalise the process described above for $n = 2$.
 
Let $z = n(x) a(y) \in \SL_n(\R)$ be a matrix in Iwasawa form.
We recall for the convenience of the reader that
\begin{displaymath}
	a(y) = \diag(d_1, \ldots, d_n) = d \cdot \begin{pmatrix}
		y_1 \cdots y_{n-1} & &  & \\
		& \ddots &  & \\
		& & y_1 & \\
		& & & 1
	\end{pmatrix}.
\end{displaymath} 
Let $\gamma \in \mathcal{M}_n(\Z, N)$ with $\det \gamma = m$ and
\begin{displaymath}
z^{-1} \gamma z = O(m^{1/n}).
\end{displaymath}
We can now multiply the previous equation with its transpose and obtain
\begin{equation}\label{eq-cond-gram}
z^{-1} \cdot \gamma \cdot z \cdot z^{T} \cdot \gamma^T \cdot z^{-T} = O(m^{2/n}).
\end{equation}
Notice now that $\gamma \cdot z \cdot z^T \cdot \gamma^T$ is the Gram matrix of the rows of $\gamma$ with respect to the scalar product defined by $z$.

Denote the rows of $\gamma$ by $\gamma_1, \ldots, \gamma_n$, and denote the rows of $n(x)^{-1} \gamma$ by $v_1, \ldots, v_n$. We compute that
\begin{displaymath}
z^{-1} \gamma z z^{T} \gamma^T z^{-T} = 
\begin{pmatrix}
\norm{v_1}^2_z \cdot d_1^{-2} & \langle v_1, v_2 \rangle_z \cdot (d_1 d_2)^{-1} & \ldots & \langle v_1, v_n\rangle_z \cdot (d_1 d_n)^{-1} \\
\ast & \norm{v_2}^2_z \cdot d_2^{-2} & \ldots & \langle v_2, v_n\rangle_z \cdot (d_2 d_n)^{-1} \\
\vdots & \vdots & \ddots & \vdots \\
\ast & \ast & \ldots  & \norm{v_n}^2_z \cdot d_n^{-2}
\end{pmatrix},
\end{displaymath}
where the matrix should be completed by noting that it is symmetric. Observe now that the condition \eqref{eq-cond-gram} reduces to
\begin{equation} \label{eq-counting-prob}
\norm{v_i}_z \ll m^{1/n} \cdot d_i,
\end{equation}
for all $i = 1, \ldots, n$, since the off-diagonal conditions simply follow by the Cauchy-Schwarz inequality. 

The strategy for counting the number of matrices $\gamma$ is to iteratively count the number of possibilities for its rows. 
More precisely, we first count the number of possible $\gamma_n = v_n$ by a lattice point counting argument using Lemma \ref{lemma-lattice-count}.
For each such fixed possibility, we then count the number of possible $\gamma_{n-1}$ by using the condition on $v_{n-1}$ in \eqref{eq-counting-prob}. 
For this observe that
\begin{displaymath}
v_{n-1} = \gamma_{n-1} - \xi \cdot \gamma_n,
\end{displaymath}
where $\xi \in \R$ can be computed from the $x$-coordinates of $z$ (in fact, $\xi = x_{n-1, n}$). 
Thus, having fixed $\gamma_n$, the condition can be interpreted as saying that $\gamma_{n-1}$ is a lattice point inside a ball with shifted centre. 
We can use that the bounds in Lemma \ref{lemma-lattice-count} are independent of the centre of the ball.
In the results below, we ultimately choose $m$ small enough so that the ball can only contain one lattice point.

We continue this process iteratively, using that $n(x)^{-1}$ is upper triangular unipotent. 
We bound the number of $\gamma$ by multiplying together the number of possibilities for each row.
As before, we only used the inequality $\det(\gamma) \leq m$ and therefore we cannot detect, at this point, the sparseness of the sequence of determinants.
This latter issue only shows up when using the unconditional amplifier and is dealt with in the next section.

To have $z$ in a good position for applying the strategy above, we make the reduction given by Proposition \ref{prop-fricke-red} and assume \textbf{case (2.b)} in the statement.

\begin{definition}
	For any $z \in \SL_n(\R)$ define
	\begin{displaymath}
		z_N = \diag(N, \ldots, N, 1) \cdot z.
	\end{displaymath}
\end{definition}

\begin{lemma} \label{lemma-balanced}
	Let $N$ be a prime and let $z \in \mathbb{H}$ satisfy the conditions of \textbf{case (2.b)} in Proposition \ref{prop-fricke-red}.
	In particular, $d = \alpha(z) = \norm{e_n}_z$ and $z \in \mathcal{L}'(B,A)$.
	Then the $\lambda_1(L_{z_n}) = d = \norm{e_n}_{z_N}$.
\end{lemma}
\begin{proof}
	Since $z \in \mathcal{L}'(B,A)$, we have $\lambda_1(L_{z'}^\ast) = d N^{-1+1/n}$.
	We now compute that
	\begin{displaymath}
		(z')^{-T} = A_N^{-1} z = N^{1/n} \cdot \diag(1, \ldots, 1, N)^{-1} z = N^{-1 + 1/n} z_N.
	\end{displaymath}
\end{proof}

The following is the main and simplest counting result of this paper and implements the strategy discussed above.

\begin{proposition} \label{prop-counting}
	Let $N$ be a prime and let $z \in \mathbb{H}$ satisfy the conditions of \textbf{case (2.b)} in Proposition \ref{prop-fricke-red}.
	Then 
	\begin{displaymath}
		|\{\gamma \in \mathcal{M}_n(\Z, N) \mid \det(\gamma) \ll \Lambda^n,\, z^{-1} \gamma z = O(\Lambda) \}|
		\ll_{n, \Omega} \Lambda^n (1 + \Lambda^n/\sqrt{N})^{n-1}.
	\end{displaymath}
\end{proposition}
\begin{proof}
	The bottom row $e_n \cdot \gamma z$ has congruence conditions and thus lies in the lattice corresponding to $z_N$. 
	By Lemma \ref{lemma-balanced}, the minimum of $L_{z_N}$ is equal $\alpha(z) = d$ and $e_n z_N$ is a vector of shortest length in $L_{z_N}$.
	
	Recall now the condition
	\begin{displaymath}
		\norm{\gamma_n}_{z} \ll \Lambda d_n = \Lambda d
	\end{displaymath}
	from \eqref{eq-counting-prob}. By Lemma \ref{lemma-lattice-count}, there are at most
	\begin{displaymath}
		1 + \frac{\Lambda d}{d} + \frac{(\Lambda d)^2}{d^2} + \cdots + \frac{(\Lambda d)^n}{d^n} \ll_{n} \Lambda^n
	\end{displaymath}
	possibilities for the row $\gamma_n = e_n \cdot \gamma$, using the trivial bounds $\lambda_i \geq \lambda_1$ for all $i$ for successive minima.
	
	We continue bounding the number of possibilities for $\gamma_i$ inductively, $i < n$. 
	More precisely, we suppose that $\gamma_{j}$ with $i < j \leq n$ are fixed. 
	Then, by using the fact that $n(x)^{-1}$ is unipotent upper triangular in condition \eqref{eq-counting-prob}, the number of possibilities left for $\gamma_i$ is bounded by the number of lattice points in $L_z$ in a ball of radius $L \cdot d_i$ with fixed centre determined by the $\gamma_{j}$, $i < j$, and $n(x)$.
	
	Next, note that $L_z$ is balanced, meaning that the successive minima of $L_z$ are all $\asymp_{n, \Omega} 1$, as in Lemma \ref{lemma-red-to-cpt}. 
	Furthermore, 
	\begin{displaymath}
		d_i^n = (d y_1 \cdots y_{n-i})^n \ll_{n, \Omega} 1/\sqrt{N},
	\end{displaymath}
	by assumption.
	By Lemma \ref{lemma-lattice-count}, there are at most 
	\begin{displaymath}
		\ll_{n, \Omega} 1 + \Lambda d_i + \cdots + (\Lambda d_i)^n \ll_{n, \Omega} 1 + \frac{\Lambda ^n}{\sqrt{N}}
	\end{displaymath}
	possibilities for $\gamma_i$. 
	
	Putting all bounds together, we bound the number of matrices $\gamma$ by
	\begin{displaymath}
		\ll_{n, \Omega} \Lambda^n (1 + \Lambda^n / \sqrt{N})^{n-1}.
	\end{displaymath}
\end{proof}

\begin{remark} \label{rem-last-row}
	The last part of the proof above shows that, as long as $\Lambda$ is small enough in terms of $N$, the choice of last row of $\gamma$ already determines the whole matrix.\footnote{
		We also remark that numerical experiments in dimension $n = 2$ seem to indicate that the bound we obtain for the possibilities for the last row might be sharp.
	}
\end{remark}

\subsection{Detecting determinants that are higher powers}\label{sec-power-dets}

The bound supplied by Proposition \ref{prop-counting} is too weak to suffice in the unconditional amplifier, Proposition \ref{prop-amp-uncond}, where powers $\nu > 1$ show up and introduce sparseness into the average. 
Taking Remark \ref{rem-last-row} into consideration, we see that the approach in the previous section is over-counting the possibilities for the last row of $\gamma$. 
Motivated by this observation, we refine the argument by counting the lattice points $\gamma_n$ only up to the contribution of the vector $e_n$. 
This latter contribution and the shape of the determinant (being a $\nu$-th power) give rise to a diophantine equation that has the right amount of solutions in the generic case. 
We then consider the degenerate case separately. 
To simplify the latter, we eventually make the assumption that the degree $n$ is prime.

For talking about the non-degenerate case, denote by $\chi_\gamma(X) = \det(X \cdot \id_n - \gamma)$ the characteristic polynomial of a matrix $\gamma$. 
We call $\gamma \in M_n(\Q)$ \emph{non-degenerate} if the polynomials 
\begin{displaymath}
	(-1)^n \chi_\gamma(X) - Y^\nu \in \Q[X, Y]
\end{displaymath} 
are irreducible over $\Q$ for all $1 \leq \nu \leq n$. Define
\begin{displaymath}
	H_\ast (z, m, N) = \{ \gamma \in H(z, m, N) \mid \gamma \text{ non-degenerate} \}.
\end{displaymath}

\begin{proposition} \label{prop-non-deg}
	Let $N$ be a prime and let $z \in \mathbb{H}$ satisfy the conditions of \textbf{case (2.b)} in Proposition \ref{prop-fricke-red}. 
	Additionally, let $L \ll N^{1/2n^2 - \varepsilon}$ and $N \gg_\Omega 1$ be large enough. Then 
	\begin{displaymath}
		\sum_{m \asymp L^n} |H_\ast(z, m^{\nu}, N)| \ll L^{(n-1)\nu} \cdot L^{1 + \varepsilon}
	\end{displaymath}
	for any $1 \leq \nu \leq n$.
\end{proposition}
\begin{proof}
	Let $\gamma \in H_\ast(z, m^{\nu}, N)$, $m \asymp L^n$, and consider again the number of possibilities for the last row $\gamma_n$. 
	For this, let $b_1, \ldots, b_{n-1}, e_n z$ be a reduced basis for $L_{z_N}$ (see Section \ref{sec-Iwasawa} for the definition, which we apply to $N^{-(n-1)/n}z_N \in \SL_n(\R)$). 
	From Lemma \ref{lemma-balanced}, we note again that $e_n$ is a vector of shortest length in $L_{z_N}$, where $\norm{e_n}_{z_N} = d$.
	We have $\norm{b_i} \geq d$.
	
	Now $\gamma_n z \in L_{z_N}$, so it can be written as 
	\begin{displaymath}
		\gamma_n z = \sum_{i=1}^{n-1} a_i b_i + a_n e_n z_N
	\end{displaymath}
	with $a_i \in \Z$.
	By Lemma \ref{lemma-red-basis-orthog} and recalling the condition $\norm{\gamma_n}_z \ll L^\nu d$ from \eqref{eq-counting-prob}, we deduce that $a_i \ll L^\nu$, for all $ 1 \leq i \leq n$.
	
	There are $L^{(n-1) \nu}$ possibilities for $a_1 ,\ldots, a_{n-1}$.
	Choose any such combination of coefficients and assume there exist $\gamma \in H_\ast(z, m^{\nu}, N)$ and $\gamma' \in H_\ast(z, l^{\nu}, N)$, for $m, l \asymp L^n$, such that $\gamma_n = \sum_{i=1}^{n-1} a_i b_i + a_n e_n z$ and 
	\begin{displaymath}
		\gamma'_n - \gamma_n = \lambda e_n.
	\end{displaymath}
	Then $\lambda \in \Z$ and $\lambda \ll L^\nu$. 
	Observe also that the matrix $\gamma - \gamma'$ satisfies the same geometric conditions \eqref{eq-counting-prob} as $\gamma$ and $\gamma'$, simply by the triangle inequality (with a doubled implied constant, of course).
	
	We now apply the same iterative process as in the proof of Proposition \ref{prop-counting}. We note however that, under the present conditions, each step yields at most one possibility. 
	Indeed, fix the last row of $\gamma - \gamma'$, having the form $\lambda e_n$, by fixing $\lambda \ll L^\nu$. 
	Next, the number of possibilities for the row $(\gamma -\gamma')_{n-1}$ is bounded by the number of $L_z$-lattice points in a ball of radius $L^\nu \cdot N^{-1/2n}$ centred at $x_{n-1, n} \cdot \lambda e_n z$, where $x_{n-1, n}$ is one of the $x$-coordinates of $z$. 
	By assumption, the radius is bounded by $N^{-\varepsilon}$. 
	However, if $N$ is large enough, this is greater than the first successive minimum of $z$, which is $\asymp_\Omega 1$. 
	There is thus only one possible lattice point.
	
	On the other hand, it is clear that the multiple $\lambda \cdot \id_n$ of the identity matrix lies in the set $H(z, \lambda^n, N)$. 
	Since $\lambda \ll L^\nu$, we see that $\lambda e_{n-1} \cdot z$ satisfies the condition of the lattice point above (again, condition \eqref{eq-counting-prob}). 
	Consequently, it follows that 
	\begin{displaymath}
		(\gamma - \gamma')_{n-1} = \lambda \cdot e_{n-1}.
	\end{displaymath}
	Iterating this argument and keeping in mind the computations in the proof of Proposition \ref{prop-counting}, we deduce that
	\begin{displaymath}
		\gamma - \gamma' = \lambda \cdot \id_n.
	\end{displaymath}
	
	It remains to count the possibilities for $\lambda$. Considering the determinant of $\gamma'$, we have
	\begin{displaymath}
		l^\nu = \det(\gamma') = \det(\gamma - \lambda \cdot \id_n) = (-1)^n \chi_\gamma(\lambda).
	\end{displaymath}
	Therefore, $(\lambda, l) \in \Z$ are a solution to the equation 
	\begin{displaymath}
		(-1)^n \chi_\gamma(X) - Y^\nu = 0.
	\end{displaymath}
	Since this polynomial is defined over $\Z$ and irreducible over $\Q$ by assumption, we count the number of such solutions using Heath-Brown's Theorem 3 in \cite{HB}. 
	In the notation there, after homogenising the polynomial, we set $B_1 = L^\nu$ for the bound on $\lambda$, then $B_2 = L^n$ for the bound on $l$, and finally $B_3 = 1$ for the bound on the additional variable. 
	Then we compute $T = L^{n \nu}$ and $V = L^{\nu + n}$. 
	Heath-Brown's result then gives the bound
	\begin{displaymath}
		\frac{V^{1/n + \varepsilon}}{T^{1/n^2}}  = L^{1 + \varepsilon}
	\end{displaymath}
	on the number of solutions we are considering. This bounds in particular the number of possibilities for $\lambda$ over all relevant determinants and so finishes the proof.
\end{proof}

We are now left with counting degenerate matrices. 
This is reminiscent of treating the special case of parabolic matrices in \cite[Lemma 2]{HT2}. 
For this we restrict to prime degrees, allowing for a clean classification of the degenerate case.

Let $n \geq 2$ be prime. Since $\chi_\gamma$ is a polynomial of degree $n$ over $\Q$, a result of Schinzel \cite{schinzel} shows that 
\begin{displaymath}
	(-1)^n \chi_\gamma(X) - Y^\nu
\end{displaymath}
is irreducible, unless $\nu = n$ and
\begin{displaymath}
	\chi_\gamma(X) = \alpha (X-\beta)^n
\end{displaymath}
for $\alpha, \beta \in \Q$. 
In the first case, it is irreducible over $\CC$ if and only if it is irreducible over $\Q$. 
In the latter case, we have $\alpha = 1$ by normalisation and $\beta^n = \det(\gamma)$.

The irreducibility criterion above and Proposition \ref{prop-non-deg}, by following its proof again verbatim, imply the following bounds.

\begin{corollary} \label{cor-count}
	Assume the same conditions as in Proposition \ref{prop-non-deg} and, additionally, let $n$ be prime.
	Then 
	\begin{displaymath}
	\sum_{m \asymp L^n} |H(z, m^{\nu}, N)| \ll L^{(n-1)\nu} \cdot L^{1 + \varepsilon}
	\end{displaymath}
	for any $1 \leq \nu \leq n-1$.
\end{corollary}

We have thus reduced the problem to counting matrices $\gamma \in H(z, m^n, N)$ for some $m \asymp L^n$, such that
\begin{displaymath}
	\chi_\gamma(X) = (X - \beta)^n.
\end{displaymath}
Since $\beta \in \Q$, it follows that $\beta = \pm m \in \Z$ (there is no sign for odd $n$).
Denote the subset of such matrices by $H_{\text{par}}(z, m^n, N)$.

The method of proof in Proposition \ref{prop-non-deg} provides even more. 
We recall at this point that the determinants $m^\nu$ appearing in the counting problem have a particular shape, namely $m = p \cdot q^{n-1}$, where $p$ and $q$ are primes of size $L$ (see the amplifier in Proposition \ref{prop-amp-uncond}).
We are thus averaging over a set of size $L^2$.
However, we can consider the special case $p = q$, which gives a smaller number of terms.

\begin{corollary} \label{prop-degenerate-p=q}
	Assume the same conditions as in Corollary \ref{cor-count}. Then
	\begin{displaymath}
		\sum_{p \asymp L} |H_{\text{par}}(z, p^{n^2}, N)| \ll L^{(n-1)n} \cdot L.
	\end{displaymath}
\end{corollary}
\begin{proof}
	We follow the proof of Proposition \ref{prop-non-deg}, but first we fix the determinant $p^{n^2}$, where $p \asymp L$. 
	There are, of course, at most $L$ such determinants. 
	Now the number of choices for a potential last row of $\gamma \in H_{\text{par}}(z, p^{n^2}, N)$ up to the contribution of $e_n$, i.e.\ up to the last component, is bounded by $L^{(n-1)n}$.
	Choose $\gamma$ and $\gamma'$ two matrices in $H_{\text{par}}(z, p^{n^2}, N)$ with the same last row up to the last component.
	
	As in the proof of Proposition \ref{prop-non-deg}, we find that
	\begin{displaymath}
		\gamma - \gamma' = \lambda \cdot \id_n.
	\end{displaymath}
	We apply again the determinant to this equation and obtain that
	\begin{displaymath}
		(\lambda - p^n)^n = p^{n^2}.
	\end{displaymath}
	It follows that there are only two possibilities for $\lambda$ and this proves the statement.
\end{proof}

We observe that the actual average of size $L^2$ would have given a bound of the form $L^{n(n-1)} \cdot L^2$, which is on the edge of what is needed for a saving. The next section significantly refines the argument to treat this issue.

\subsection{Counting at different cusps}

Corollary \ref{prop-degenerate-p=q} allows us now to reduce the problem further. We are now counting matrices $\gamma$ in the set
\begin{displaymath}
	\bigcup_{\substack{p, q \asymp L \\ p \neq q}} H_{\text{par}}(z, (pq^{n-1})^n, N).
\end{displaymath} 
By Theorem III.12 in \cite{newman}, there is $h \in \SL_n(\Z)$ such that
\begin{equation} \label{eq-upper-tri}
	h\gamma h^{-1} = \begin{pmatrix}
	m & \ast &  \ast\\
	  & \ddots    & \ast\\
	  &      & m
	\end{pmatrix}
\end{equation}
is upper triangular with $m$ on the diagonal. Indeed, $\chi_\gamma$ splits into linear factors and thus the blocks in \cite[Thm. III.12]{newman} are one dimensional.

In the simplest case, we could assume that $h \in \Gamma_0(N)$. 
The next lemma shows that this is almost the same as assuming that $h = 1$ and that $\gamma$ has the same last row as the identity matrix, in which case we apply the philosophy from Remark \ref{rem-last-row}, namely that the last row determines the matrix. 
However, we remark here already that there are other possibilities for $h$ that correspond to different \emph{cusps}, as in Lemma \ref{lemma-cusps} below, for which counting becomes more difficult.

\begin{lemma} \label{lemma-upper-tri-id}
	Assume the same conditions as in Proposition \ref{prop-non-deg} and let $\gamma \in H(z, m^n, N)$ or $\gamma \in H(z', m^n, N)$ for $m \asymp L^n$. 
	If there exists $h \in \Gamma_0(N)$ such $h \gamma h^{-1}$ has last row equal to $m \cdot e_n = (0, \ldots, 0, m)$, then $\gamma = m \cdot \id_n$.
\end{lemma}
\begin{proof}
	Assume that $\gamma \in H(z, m^n, N)$. Since $h \in \Gamma_0(N)$, it is easy to see from the definition that $\gamma \in H(z, m^n, N)$ implies $\eta := h \gamma h^{-1} \in H(hz, m^n, N)$. 
	Consider the Iwasawa coordinates of $hz = n(x)a(y)$. 
	Multiplying $h$ from the left by a matrix of the form
	\begin{displaymath}
		\begin{pmatrix}
		\xi & \\
		& 1
		\end{pmatrix} \in \Gamma_0(N)
	\end{displaymath}
	with $\xi \in \SL_n(\Z)$, we may assume that $y_i \gg 1$ for $i = 2, \ldots, n-1$ (see Remark \ref{rem-reduction-n-1}). 
	Under such a modification, we may also still assume that the last row $e_n \eta$ has the form $(0, \ldots, 0, m) = m\cdot e_n$.

	To obtain from this bounds on the entries of $a(y)$ we note that, since $\det(hz) = 1$,
	\begin{displaymath}
		\norm{e_2 \wedge \cdots \wedge e_n}_{h z} = (d y_1 \cdots y_{n-1})^{-1} \geq \beta(z),
	\end{displaymath}
	recalling Remark \ref{rem-iwasawa-calculations}, the definition of $\beta(z)$ in Section \ref{sec-fricke-red}, and that $h \in \Gamma_0(N)$. 
	We are assuming that $z$ satisfies the conditions in case \textbf{(2.b)} of Proposition \ref{sec-fricke-red}.
	Thus, $z \in \mathcal{L}'(B, A)$, which by Table \ref{table} implies that $\beta(z) = N^{1/n} \lambda_1(L_{z'})$.
	The assumptions also give that $\lambda_1(L_{z'}) \gg N^{-1/2n}$ and so $\beta(z) \gg N^{1/2n}$. 
	As such, we have
	\begin{displaymath}
		dy_1 \ll dy_1 y_2 \ll \ldots \ll d y_1 \cdots y_{n-1} \ll N^{-1/2n}.
	\end{displaymath}
	
	We are now again in a similar situation as in the counting results Proposition \ref{prop-counting} and Proposition \ref{prop-non-deg}. 
	However, the last row of $\eta$ is already fixed to be $m \cdot e_n$. 
	As in Proposition \ref{prop-non-deg}, the assumption $L \ll N^{1/2n^2 - \varepsilon}$, the bound above on the entries of $a(y)$, and the assumption that $L_z$ is balanced imply that the last row of $\eta$ determines the whole matrix. 
	Therefore, $\eta = m \cdot \id_n$ and so, undoing conjugation, $\gamma = m \cdot \id_n$.
	
	The case $\gamma \in H(z', m^n, N)$ follows analogously. 
	Here $z \in \mathcal{L}(B, A)$, but $L_z$ is also balanced, so $\beta(z') = N^{1/n} \lambda_1(L_{z}) \asymp N^{1/n}$.
	Therefore, the $y$-coordinates of $hz'$ satisfy the stronger bound
	\begin{displaymath}
		dy_1 \ll dy_1 y_2 \ll \ldots \ll d y_1 \cdots y_{n-1} \ll N^{-1/n}.
	\end{displaymath}
	On the other hand, $L_{z'}$ is potentially unbalanced, yet $\lambda_1(L_{z'}) \gg N^{-1/2n}$.
	By following the proof of Proposition \ref{prop-non-deg}, we are considering $z'$-lattice points inside balls of radius $L^{n} \cdot d y_1 \cdots y_k$ for some $k$.
	The assumption on $L$ bounds the radius by $N^{1/2n -\varepsilon - 1/n} = N^{-1/2n - \varepsilon}$, which is smaller than the shortest length in $L_{z'}$ for $N$ large enough.
	Thus, there is at most on lattice point in such a ball, and the proof continues as usual.
\end{proof}

We investigate now the cusps of $\Gamma_0(N)$ with respect to the minimal parabolic. Define therefore $U_n(\Z)$ to be the subgroup of $\SL_n(\Z)$ of unipotent upper triangular matrices, that is, with ones on the diagonal. 

Let also $W_n \leq \SL_n(\Z)$ denote the subgroup of permutation matrices. We call two such matrices equivalent if they have the same last row and denote by $\overline{W}_n$ the set of equivalence classes. By considering $\SL_{n-1}(\Z)$ embedded inside $\Gamma_0(N)$, it is easy to see that
\begin{displaymath}
	\overline{W}_n \cong (\Gamma_0(N) \cap W_n) \backslash W_n.
\end{displaymath}
Note also that $|\overline{W}_n| = n$.

\begin{lemma} \label{lemma-cusps}
	Let $N$ be prime. 
	Then any system of representatives for $\overline{W}_n$ is a system of representatives for the double quotient
	\begin{displaymath}
		\Gamma_0(N) \backslash \SL_n(\Z) / U_n(\Z).
	\end{displaymath}
\end{lemma}
\begin{proof}
	Let $\xi \in \SL_n(\Z)$ and let $(a_1, \ldots, a_n)$ be the first column of $\xi$, a primitive vector in $\Z^n$. 
	First, we reduce $a_n$ to either $0$ or $1$ by acting from the left by $\Gamma_0(N)$.
	
	Indeed, assume that $\gcd(a_n, N) = 1$. 
	Then the vector $(Na_1, \ldots, Na_{n-1}, a_n)$ is also primitive. 
	Therefore, there is a primitive $(b_1, \ldots, b_n) \in \Z^n$ such that
	\begin{displaymath}
		N a_1 b_1 + \ldots + N a_{n-1} b_{n-1} + a_n b_n = 1.
	\end{displaymath}
	From this it is clear that $\gcd(N, b_n) = 1$ so that $(N b_1, \ldots, N b_{n-1}, b_n)$ is primitive. 
	Let $\gamma \in \SL_n(\Z)$ be a matrix with the latter as its last row. 
	Then $\gamma \in \Gamma_0(N)$ and $\gamma \xi$ has last row of the form $(1, \ast, \ldots, \ast)$.
	
	Since $N$ is prime, negating the assumption above means that $N \mid a_n$. 
	Now let $d = \gcd(a_1, \ldots, a_{n-1})$. 
	Then $\gcd(a_n, d) = 1$ and there exists a primitive vector $(b_1, \ldots, b_{n-1})$ such that
	\begin{displaymath}
		b_1 a_1 + \ldots b_{n-1} a_{n-1} = d.
	\end{displaymath}
	Therefore,
	\begin{displaymath}
		\sum_{i = 1}^{n-1} (a_n b_i) \cdot a_i + (-d) a_n = 0.
	\end{displaymath}
	The vector $(a_n b_1, \ldots, a_n b_{n-1}, -d)$ is primitive by the observations above, so there is $\gamma \in \SL_n(\Z)$ with this vector as its last row. 
	Again, $\gamma \in \Gamma_0(N)$ since $N \mid a_n$ and the last row of $\gamma \xi$ has the form $(0, \ast, \ldots, \ast)$.
	
	Assume now that $\xi$ has last row of the form $(1, \ast, \ldots, \ast)$. 
	It is clear that we can multiply $\xi$ from the right by a matrix in $U_n(\Z)$ such that the resulting last row is simply $(1, 0, \ldots, 0)$.
	Call this new matrix $\xi$ again and take $w \in W_n$ a permutation matrix with the same last row (for instance the so-called long Weyl element).
	In other words, $e_n \xi = e_n w$, where $e_n$ is the $n$-th standard basis vector $(0, \ldots, 0, 1)$.
	The matrix $w \xi^{-1}$ preserves $e_n$ so it must have $e_n$ as its last row. 
	In particular, $w \xi^{-1} \in \Gamma_0(N)$ and we are done in this case.
	
	On the other hand, let $\xi$ have last row of the form $(0, \ast, \ldots, \ast)$. 
	Using the embedding of $\SL_{n-1}(\Z)$ in the upper left corner of $\Gamma_0(N)$, we may modify $\xi$ so that its first column is of the form $(1, 0, \ldots, 0)$, by similar arguments. 
	This now allows an inductive procedure, considering the lower right $(n-1) \times (n-1)$ block of $\xi$.
	We see that one can always reduce the last row of $\xi$ to be a standard basis vector and the paragraph above shows how to obtain a permutation matrix from $\xi$.
	
	To check that no two such representatives in $\overline{W}_n$ produce the same double coset is easy. 
	For $w_1, w_2 \in W_n$, if $w_1 = \gamma w_2 u$ with $\gamma \in \Gamma_0(N)$ and $u \in U_n(\Z)$, then $\gamma = w_1 u^{-1} w_2$.
	One now computes the shape of $U_n(\Z)$ transformed by permutation of rows and of columns. 
	We leave out the details of this argument.
\end{proof}

\begin{remark} \label{rem-weyl-el}
	We make the following simple observation that becomes very useful in the arguments below. Let $w_k \in \overline{W}_n$ be a representative with last row equal to $e_k$. We can take $w_n = \id_n$. We can also take $w_1$ to be the long Weyl element
	\begin{displaymath}
		w_1 = \begin{pmatrix}
			& & & & 1 \\
			& & & 1 & \\
			& & \reflectbox{$\ddots$} & & \\
			& 1 & & & \\
			1 & & & &
		\end{pmatrix}
	\end{displaymath}
	with ones on the anti-diagonal. Finally, for any $k \neq 1$, we can choose the representative $w_k$ to have first row (and thus also first column) equal to $e_1$.
\end{remark}

We finally state the main result for degenerate matrices below and recall the additional condition on the determinantal divisors appearing in the amplifier, Proposition \ref{prop-amp-uncond}.

\begin{proposition} \label{prop-degenerate-count}
	Assume the same conditions as in Corollary \ref{cor-count}. For $N$ large enough, the set of matrices $\gamma$ possibly occurring in $H_{\text{par}}(z, (pq^{n-1})^n, N)$ for some primes $p,q \asymp L$, $p \neq q$, satisfying additionally that
	\begin{displaymath}
		\Delta_{n-1}(\gamma) = q^{n(n-2)}
	\end{displaymath} 
	is empty.
\end{proposition}

It is perhaps useful at this point to give a brief overview of the proof. 
We make a case distinction, based on the cusp classification above. 
If $h$ in \eqref{eq-upper-tri} corresponds to the identity $w_n$, then we are done by Lemma \ref{lemma-upper-tri-id}. 
If $h$ corresponds to the long Weyl element $w_1$, we apply the Fricke involution, which effectively switches the cases $w_n$ and $w_1$, and so the same lemma, available for both $z$ and $z'$, finishes this case as well.

In fact, the Fricke involution generally exchanges the cases $w_k$ and $w_{n+1-k}$.
\begin{displaymath}
	\begin{tikzcd}
		& {w_2} & \cdots & {w_{n-1}} \\
		{w_1} &&&& {w_n}
		\arrow["{W_N}", curve={height=12pt}, tail reversed, from=2-1, to=2-5]
		\arrow[curve={height=12pt}, tail reversed, from=1-2, to=1-4]
	\end{tikzcd}
\end{displaymath}
However, the usual counting argument, choosing vectors step-by-step from the bottom of the matrix going upwards, seems difficult to implement in the intermediate cases $1 < k < n$.
It is here that the assumption $p \neq q$, together with the seemingly harmless choice of representatives $w_k$ in Remark \ref{rem-weyl-el}, comes in.
Indeed, the choice of representatives is akin to a very weak balancedness assumption on the new, unknown basis for the lattice that appears in the counting problem.
This assumption implies that at least one element of the superdiagonal of the upper triangular matrix in \eqref{eq-upper-tri} is zero.
Computing $\Delta_{n-1}$, this is enough to derive a contradiction to $p \neq q$.

\begin{proof}
	Let $m = p q^{n-1}$. As in \eqref{eq-upper-tri}, there is $h \in \SL_n(\Z)$ such that $h \gamma h^{-1}$ is upper-triangular with diagonal $(m, \ldots, m)$. By Lemma \ref{lemma-cusps} we can write $h^{-1} = \sigma^{-1} w u^{-1}$ with $\sigma \in \Gamma_0(N)$, $u \in U_n(\Z)$, and $w \in \overline{W}_n$.
	
	Next, conjugating by $u$, we easily see that
	\begin{equation} \label{eq-eta}
		w^T \sigma \gamma \sigma^{-1} w = \begin{pmatrix}
			m & \ast &  \ast\\
			& \ddots    & \ast\\
			&      & m
		\end{pmatrix} =: \eta
	\end{equation}
	is also of the same form. 
	Now if $w = w_n = \id_n$, meaning that the last row of $w$ is $e_n$ as in Remark \ref{rem-weyl-el}, we are done by Lemma \ref{lemma-upper-tri-id}. 
	The latter implies that $\gamma = pq^{n-1} \cdot \id_n$, which does not have the required determinantal divisors and leads to a contradiction.
	
	If $w = w_1$ is the long Weyl element, we apply the Fricke involution. By transposing the condition 	
	\begin{displaymath}
		z^{-1} \gamma z = O(m),
	\end{displaymath}
	we see that
	\begin{displaymath}
		A_N (\sigma \gamma \sigma^{-1})^T A_N^{-1}
	\end{displaymath}
	lies in $H(\tilde{\sigma} z', m^n, N)$ with some $\tilde{\sigma} \in \Gamma_0(N)$.
	
	Next, observe that
	\begin{displaymath}
		(\sigma \gamma \sigma^{-1})^T = w \eta^T w^T
	\end{displaymath}
	is again upper triangular. 
	By Lemma \ref{lemma-upper-tri-id}, we deduce that
	\begin{displaymath}
		A_N (\sigma \gamma \sigma^{-1})^T A_N^{-1} = m \cdot \id_n
	\end{displaymath}
	and thus $\gamma = m \cdot \id_n$, which is a contradiction again.
	
	Finally, let $w = w_k$ with $1 < k < n$. 
	Notice first that \eqref{eq-eta} and the congruences modulo $N$ satisfied by $\gamma$ and $\sigma$ imply that the $k$-th row of $\eta$ also satisfies congruences.
	Indeed, $w_k \eta w_k^T$ is a matrix of the $\Gamma_0(N)$ shape.
	More precisely, $N \mid \eta_{kj}$ for $j > k$. Since $k < n$, we have in particular $N \mid \eta_{k, k+1}$.
	
	Let us now assume that the superdiagonal of $\eta$ only contains non-zero elements. 
	That is, $\eta_{j, j+1} \neq 0$ for all $1 \leq j < n$. Recall the condition
	\begin{displaymath}
		z^{-1} \gamma z = (w^T \sigma z)^{-1} \eta (w^T \sigma z) = O(m).
	\end{displaymath}
	We can rewrite $w^T \sigma z = n(x) \cdot a(y)$ in Iwasawa coordinates (indeed, conjugating by an orthogonal matrix leaves $O(m)$ invariant), denoting the $y$-coordinates as usual. 
	It is now a common and important observation that the superdiagonal of upper triangular matrices enjoys a certain additive abelian-like property with respect to matrix multiplication. 
	This observation or direct computation should convince the reader that
	\begin{displaymath}
		(na)^{-1} \eta na = \begin{pmatrix}
			m & y_{n-1}^{-1} \eta_{1,2} & \ast &  \ldots \\
			& m & y_{n-2}^{-1} \eta_{2,3}    & \ldots \\
			&  & \vdots & \vdots \\
			& & m & y_1^{-1} \eta_{n-1, n} \\
			& & & m
		\end{pmatrix}.
	\end{displaymath}
	Since this is $O(m)$, the assumption that $|\eta_{j, j+1}| \geq 1$ now implies that $y_j \gg 1/m$. 
	Even more and crucially, recall that $N \mid \eta_{k, k+1}$, so that $y_{n-k} \gg N/m$. 
	Putting these together, we obtain the bound
	\begin{displaymath}
		y_1 \cdots y_{n-1} \gg \frac{N}{m^{n-1}} \gg N^{1/2}
	\end{displaymath}
	using the assumption $L \ll N^{1/2n^2 - \varepsilon}$ and that $m \asymp L^n$. 
	
	We return now to a technique used in the proof of Lemma \ref{lemma-upper-tri-id}. 
	We observe again that
	\begin{displaymath}
		\norm{e_2 \wedge \cdots \wedge e_n}_{w^T \sigma z} = (d y_1 \cdots y_{n-1})^{-1}.
	\end{displaymath}
	On the other hand, our choice of representative $w = w_k$ in Remark \ref{rem-weyl-el} implies that the first row of $w^T$ is equal to $e_1$ and the other rows are permuted between them in some way.
	This means that
	\begin{displaymath}
		\norm{e_2 \wedge \cdots \wedge e_n}_{w^T \sigma z} = \norm{e_2 \wedge \cdots \wedge e_n}_{\sigma z} \geq \beta(z).
	\end{displaymath}
	Recall that, in case \textbf{(2.b)} of Proposition \ref{prop-fricke-red}, we have $\beta(z) = N^{1/n} \lambda_1(L_{z'})$.
	The assumptions also give that $\lambda_1(L_{z'}) \gg N^{-1/2n}$, so $\beta(z) \gg N^{1/2n}$.
	We obtain that
	\begin{displaymath}
		d y_1 \cdots y_{n-1} \ll N^{-1/2n}.
	\end{displaymath}
	Recall also that $d = \norm{e_n}_{w^T \sigma z}$, and since $z$ defines a balanced lattice, $d \gg 1$. 
	Therefore,
	\begin{displaymath}
		y_1 \cdots y_{n-1} \ll N^{1/n},
	\end{displaymath}
	which constitutes a contradiction to the previous paragraph for large enough $N$.
	
	We deduce that the superdiagonal of $\eta$ must contain some zero. 
	It is now straight-forward to prove that $m$ divides $\Delta_{n-1}(\eta)$. 
	Indeed, the only $(n-1) \times (n-1)$ minor that is not obviously divisible by $m$ is the upper right minor, formed by removing the first column and the last row of $\eta$. 
	Proving the claim here is an easy exercise in Laplace, or cofactor, expansion. 
	
	Observe now that the invariance properties of determinantal divisors (see \cite[Thm. II.8]{newman}) imply that
	\begin{displaymath}
		\Delta_{n-1}(\eta) = \Delta_{n-1}(\gamma),
	\end{displaymath}
	since $w, \sigma \in \SL_n(\Z)$. Since $p \mid m$, it follows from the paragraph above and our assumption on the determinantal divisors that
	\begin{displaymath}
		p \mid q^{n(n-2)}.
	\end{displaymath}
	If $n > 2$, this implies that $p = q$, which is a contradiction to the assumption.
\end{proof}

\begin{remark}
	Notice that the case $n=2$ does not involve any intermediate Weyl elements. 
	Indeed, there are only two cusps and both reduce as above to counting upper-triangular matrices directly. 
	A more general result (for square-free levels) is contained in a slightly different language in \cite[Lemma 4.1]{HT2}.
\end{remark}

The counting results of this section taken together produce the following corollary. 
It gives a solution to the counting problem for prime $n$ that can be successfully applied to the sup-norm problem through the amplifier in Proposition \ref{prop-amp-uncond}.

\begin{corollary} \label{cor-count-p}
	Let $n$ and $N$ be a prime, and let $z \in \mathbb{H}$ satisfy the conditions of case \textbf{(2.b)} in Proposition \ref{prop-fricke-red}.
	Let $L \ll N^{1/2n^2 - \varepsilon}$ and assume that $N \gg_{\Omega, \varepsilon} 1$ is large enough. Then
	\begin{displaymath}
	\sum_{p, q \asymp L} |H(z, p^{\nu}, q^{(n-1)\nu}, N)| \ll L^{(n-1) \nu} \cdot L^{1 + \varepsilon}
	\end{displaymath}
	for any $1 \leq \nu \leq n$.
\end{corollary}

\section{Final steps}

First, assume Hypothesis \eqref{eq-hypothesis}. Proposition \ref{prop-amp} and Proposition \ref{prop-counting} together with the prime number theorem imply that
\begin{displaymath}
	\phi(z)^2 \ll_{\mu, \Omega, \varepsilon} L^{-1/2 + \varepsilon} + L^{-1/2-n + \varepsilon} \cdot L^n(1 + L^n/N)^{n-1},
\end{displaymath}
under the assumptions on $z$ specified in Proposition \ref{prop-counting}. Optimising the size of $L$, we choose $L = N^{1/2n}$. In this case, we have
\begin{displaymath}
	\phi(z) \ll L^{-1/4 + \varepsilon} \ll N^{-1/8n + \varepsilon}.
\end{displaymath}

These bounds are valid on the subsets of $\mathbb{H}$ corresponding to case \textbf{(2.b)} in Proposition~\ref{prop-fricke-red}.
As remarked at the beginning of Section \ref{sec-red}, by applying the Fricke involution, these now extend to a subset of $\Omega_N^\natural$ that we call $\Omega_N$, the bulk. 
The complement, $C_N$, is given by the points corresponding to the cases \textbf{(1)} and \textbf{(2.a)}.
These form a set with $y$-coordinates as in \eqref{eq:exceptional}, thus away from the bulk of $\Omega_N^\natural$, of volume bounded from above by $N^{(n-1)(1-1/2n^4)}$.
We use here that the volume form is invariant under the Fricke involution.
One can see this, for instance, from the invariance of the Riemannian metric under $g \mapsto g^{-T}$.

Without assuming Hypothesis \eqref{eq-hypothesis}, we let $n$ be prime and we apply Proposition \ref{prop-amp-uncond} using the counting result Corollary \ref{cor-count-p}. Similarly to the computation above, we have
\begin{displaymath}
	\phi(z)^2 \ll_{\mu, \Omega, \varepsilon} L^{-1 + \varepsilon} + \sum_{\nu = 1}^{n} \frac{1}{L^{(n-1) \nu}} \cdot L^{(n-1) \nu} L^{1 + \varepsilon}
\end{displaymath}
for $L \ll N^{1/2n^2 - \varepsilon}$. 
Maximising $L$, we get
\begin{displaymath}
	\phi(z) \ll L^{-1/2 + \varepsilon} \ll N^{-1/4n^2 + \varepsilon}.
\end{displaymath}

\section{Appendix: The Fourier bound} \label{sec-fourier}
Consider the cases \textbf{(1)} and \textbf{(2.a)} in Proposition \ref{prop-fricke-red}.
They give lower bounds for $y$-coordinates, away from the bulk $y_1 \asymp 1/N$ and $y_i \asymp 1$ for $i \geq 2$.
In an ideal situation, as in the case of $\PGL(2)$, the Fourier bound would be effective already in such regions.
Indeed, as in \cite[Lemma~4]{HT2}, there we have $\phi(x + iy) \ll N^{-1/2+\varepsilon}y^{-1/2}$, and we obtain a sub-baseline bound as soon as $y \geq N^{-1 + \delta}$.
The limit of this tool is precisely the bulk $y \asymp 1/N$.

However, the state of the art does not provide such useful Fourier bounds.
This is due to a lack of control in the Whittaker expansion, particularly for the dual form of $\phi$, as well as the weakness of bounds on residues of Rankin-Selberg $L$-functions in higher rank.

For a demonstration of the current technology, we prove a bound in case \textbf{(1)}. 
It is unconditional, yet even under a no-Siegel-zero assumption there are serious obstacles to generalising and improving it.
This goes beyond the scope of this paper.

\begin{proposition} \label{prop-fourier}
	Let $\phi$ be an $L^2$-normalised Hecke-Maaß newform of prime level $N$ and spectral parameter $\mu$, and let $z \in \Omega$ for some compactum $\Omega \subset \mathbb{H}$. For $\varepsilon > 0$ we have
	\begin{displaymath}
		\phi(z) \ll_{\Omega, \mu \varepsilon} N^{-1/2 + 1/2n + \varepsilon}.
	\end{displaymath}
\end{proposition}
\begin{proof}
	We use the bound given in Theorem 3 of \cite{BHMn}, making the necessary adjustments from the level $1$ results to level $N$. 
	The proof is very similar, so we refer to \cite{BHMn} for more details and mostly remark on what changes need to be made.
	
	Note first that the method of proof involves the Whittaker expansion \cite[(46)]{BHMn}. 
	An automorphic form for the group $\Gamma_0(N)$ enjoys the same type of Whittaker expansion, since $\SL_{n-1}(\Z)$ embeds in the upper left $(n-1)\times(n-1)$ block of $\Gamma_0(N)$, so that one can follow the same arguments given in, for instance, \cite[Theorem 5.3.2]{goldfeld} in level $1$. 
	To follow the arguments in \cite{BHMn} further, we normalise $\phi$ arithmetically, so that the first coefficient in the expansion is $1$.
	
	Next, the bound \cite[(49)]{BHMn} for $L(1+\varepsilon, \pi \times \tilde{\pi})$ holds similarly, with an additional $N^\varepsilon$ on the right-hand side. 
	Here, we let $\pi$ be the automorphic representation generated by $\phi$.
	Finally, to account for the factor between arithmetically normalised forms and $L^2$-normalised forms, we note the display before \cite[(66)]{BHMn}. 
	More precisely, if we assume $\phi$ to be arithmetically normalised, as in \cite[(46)]{BHMn}, then standard Rankin-Selberg theory shows that
	\begin{displaymath}
		\norm{\phi}^2 \asymp_{\mu} \operatorname{vol}(\Gamma_0(N) \backslash \mathbb{H}) \cdot \operatorname{res}_{s=1} L(s, \pi \times \tilde{\pi}).
	\end{displaymath}
	By \cite[Theorem 3]{brumley}, as in the two displays after (43) in \cite[Appendix]{lapid}, we can use the lower bound
	\begin{displaymath}
		\operatorname{res}_{s=1} L(s, \pi \times \tilde{\pi}) \gg C(\pi \times \tilde{\pi})^{-1/2 + 1/2n - \varepsilon},
	\end{displaymath}
	where $C(\pi \times \tilde{\pi}) = C(\pi \times \tilde{\pi}, 0)$ is the analytic conductor of $L(s, \pi \times \tilde{\pi})$. 
	We have $C(\pi) \asymp_\mu N$ and by \cite[Thm.~B.1]{BTZ} the bound 
	\begin{displaymath}
		C(\pi \times \tilde{\pi}) \ll N^{n\cdot 1 + n\cdot 1 - 2}
	\end{displaymath} 
	holds.
	
	It is easy to compute that $\operatorname{vol}(\Gamma_0(N) \backslash \mathbb{H}) \asymp N^{n-1}$. Therefore,
	\begin{displaymath}
		\norm{\phi}^2 \gg_\mu N^{n-1} \cdot N^{(2n-2)(-1/2+1/2n) - \varepsilon} = N^{1 - 1/n - \varepsilon}.
	\end{displaymath}
	Going back to $\phi$ being $L^2$-normalised by putting together the bound above and \cite[(49)]{BHMn} with the indicated adjustments, we deduce the claim.
\end{proof}

\begin{remark} \label{rem-fourier}
	In fact, one expects that $\operatorname{res}_{s=1} L(s, \pi \times \tilde{\pi}) \gg N^{-\varepsilon}$. 
	This is proven by Hoffstein-Lockhart in the case $n=2$ and for this reason we have
	\begin{displaymath}
		\phi(z) \ll_{\Omega, \mu, \varepsilon} N^{-1/2},
	\end{displaymath}
	for $z \in \Omega$ as in \cite[Lemma 4]{HT2}, for example. 
	This is the best-possible bound on the global sup-norm in this case.
\end{remark}

\begin{remark}
	More precisely, for $n > 2$, Blomer--Harcos--Maga \cite[(63)]{BHMn} prove a bound of the form
	\begin{displaymath}
		\phi(z) \ll \left( \prod_{i=1}^{n-1} y_i^{(n-i)(n-i-1)/2} \right)^{-1} = y_1^{-(n-1)(n-2)/2} \cdots y_{n-2}^{-1} \cdot y_{n-1}^0,
	\end{displaymath}
	assuming that $y_1, \ldots, y_{n-1} \gg 1$, e.g. $z$ is in a standard fundamental domain for $\SL_n(\Z)$.
	Even if this would hold for smaller $y_1$ and assuming a Hoffstein-Lockhart-type bound, one would obtain a sub-baseline bound only for $y_1 \gg N^{-2/(n-2) + \varepsilon}$ and $y_i \gg 1$ for all $i \geq 2$.
	This is still very far from covering the region $C_N$.
	
	However, the form of the bound, given through a product of $y$-coordinates, does suggest decay as soon as at least one of the coordinates is significantly larger than in the bulk.
	This is strengthened in \cite[(6.5)]{BHMn} using the dual form of $\phi$, where \emph{all} $y$-coordinates appear with negative exponents.
	However, when $N > 1$, the latter has a more complicated Whittaker expansion, which deserves to be investigated.
\end{remark}

\section*{Acknowledgements}
I am very grateful to Valentin Blomer for suggesting the problem and for his generous and valuable support.
Part of this paper was worked out while I was a visiting student at EPFL, and I am grateful to Philippe Michel and his group for great working conditions. I also thank Edgar Assing, Farrell Brumley, Gergely Harcos, and Aurel Page for useful discussions on topics related to this paper, and the referee for improving readability and correctness.

\printbibliography

@Article{BM4,
  author     = {Blomer, Valentin and Maga, P\'{e}ter},
  title      = {The sup-norm problem for {PGL}(4)},
  journal    = {Int. Math. Res. Not. IMRN},
  year       = {2015},
  number     = {14},
  pages      = {5311--5332},
  issn       = {1073-7928},
  doi        = {10.1093/imrn/rnu100},
  fjournal   = {International Mathematics Research Notices. IMRN},
  mrclass    = {11F37},
  mrnumber   = {3384442},
  mrreviewer = {Guangshi L\"{u}},
  url        = {https://doi.org/10.1093/imrn/rnu100},
}

@Article{evertse,
  author  = {Evertse, Jan-Hendrik},
  title   = {Mahler’s Work on the Geometry of Numbers},
  journal = {Documenta Mathematica},
  year    = {2019},
  volume  = {Extra Volume Mahler Selecta},
  pages   = {29-43},
  doi     = {https//:doi.org/10.25537/dm.2019.SB-29-43},
}

@Book{goldfeld,
  title      = {Automorphic forms and {$L$}-functions for the group {${\rm GL}(n,\mathbb{R})$}},
  publisher  = {Cambridge University Press, Cambridge},
  year       = {2006},
  author     = {Goldfeld, Dorian},
  volume     = {99},
  series     = {Cambridge Studies in Advanced Mathematics},
  isbn       = {978-0-521-83771-2; 0-521-83771-5},
  note       = {With an appendix by Kevin A. Broughan},
  doi        = {10.1017/CBO9780511542923},
  mrclass    = {11F55 (11F66 11F70 11R39)},
  mrnumber   = {2254662},
  mrreviewer = {Emmanuel P. Royer},
  pages      = {xiv+493},
  url        = {https://doi.org/10.1017/CBO9780511542923},
}

@Article{mahler,
  author     = {Mahler, Kurt},
  journal    = {Proc. London Math. Soc. (3)},
  title      = {On compound convex bodies. {I}},
  year       = {1955},
  issn       = {0024-6115},
  pages      = {358--379},
  volume     = {5},
  doi        = {10.1112/plms/s3-5.3.358},
  fjournal   = {Proceedings of the London Mathematical Society. Third Series},
  mrclass    = {10.3X},
  mrnumber   = {74460},
  mrreviewer = {P. Scherk},
  url        = {https://doi.org/10.1112/plms/s3-5.3.358},
}

@Article{BMn,
  author     = {Blomer, Valentin and Maga, P\'{e}ter},
  title      = {Subconvexity for sup-norms of cusp forms on {$\rm {PGL}(n)$}},
  journal    = {Selecta Math. (N.S.)},
  year       = {2016},
  volume     = {22},
  number     = {3},
  pages      = {1269--1287},
  issn       = {1022-1824},
  doi        = {10.1007/s00029-015-0219-5},
  fjournal   = {Selecta Mathematica. New Series},
  mrclass    = {11F55 (11D75 11F72)},
  mrnumber   = {3518551},
  mrreviewer = {Martin Raum},
  url        = {https://doi.org/10.1007/s00029-015-0219-5},
}

@Misc{KNS,
  author    = {Khayutin, Ilya and Nelson, Paul D. and Steiner, Raphael S.},
  title     = {Theta functions, fourth moments of eigenforms, and the sup-norm problem II},
  year      = {2022},
  copyright = {arXiv.org perpetual, non-exclusive license},
  doi       = {10.48550/ARXIV.2207.12351},
  publisher = {arXiv},
  url       = {https://arxiv.org/abs/2207.12351},
}

@Book{borel,
  author    = {Borel, Armand},
  publisher = {American Mathematical Society, Providence, RI},
  title     = {Introduction to arithmetic groups},
  year      = {2019},
  isbn      = {978-1-4704-5231-5},
  series    = {University Lecture Series},
  volume    = {73},
  doi       = {10.1090/ulect/073},
  mrclass   = {22E40 (14Lxx 20G20 20G30)},
  mrnumber  = {3970984},
  pages     = {xii+118},
  url       = {https://doi.org/10.1090/ulect/073},
}

@Book{cassels,
  author    = {Cassels, J. W. S.},
  publisher = {Berlin: Springer},
  title     = {An introduction to the geometry of numbers.},
  year      = {1997},
  edition   = {Repr. of the 1971 ed.},
  isbn      = {3-540-61788-4},
  series    = {Class. Math.},
  fseries   = {Classics in Mathematics},
  issn      = {1431-0821},
  language  = {German},
  zbl       = {0866.11041},
  zbmath    = {971533},
}

@Article{BHM,
  author   = {Blomer, Valentin and Harcos, Gergely and Mili{\'c}evi{\'c}, Djordje},
  journal  = {Duke Mathematical Journal},
  title    = {Bounds for eigenforms on arithmetic hyperbolic 3-manifolds},
  year     = {2016},
  issn     = {0012-7094},
  number   = {4},
  pages    = {625--659},
  volume   = {165},
  doi      = {10.1215/00127094-3166952},
  language = {English},
  zbl      = {1339.11062},
  zbmath   = {6571492},
}

@Article{HT3,
  author     = {Harcos, Gergely and Templier, Nicolas},
  title      = {On the sup-norm of {M}aass cusp forms of large level. {III}},
  journal    = {Math. Ann.},
  year       = {2013},
  volume     = {356},
  number     = {1},
  pages      = {209--216},
  issn       = {0025-5831},
  doi        = {10.1007/s00208-012-0844-7},
  fjournal   = {Mathematische Annalen},
  mrclass    = {11F12 (11D45 14G35)},
  mrnumber   = {3038127},
  mrreviewer = {Henri Darmon},
  url        = {https://doi.org/10.1007/s00208-012-0844-7},
}

@Article{HT2,
  author     = {Harcos, Gergely and Templier, Nicolas},
  title      = {On the sup-norm of {M}aass cusp forms of large level: {II}},
  journal    = {Int. Math. Res. Not. IMRN},
  year       = {2012},
  number     = {20},
  pages      = {4764--4774},
  issn       = {1073-7928},
  doi        = {10.1093/imrn/rnr202},
  fjournal   = {International Mathematics Research Notices. IMRN},
  mrclass    = {11F12 (11D45 11G35)},
  mrnumber   = {2989618},
  mrreviewer = {Rainer Schulze-Pillot},
  url        = {https://doi.org/10.1093/imrn/rnr202},
}

@Article{BHolo,
  author     = {Blomer, Valentin and Holowinsky, Roman},
  title      = {Bounding sup-norms of cusp forms of large level},
  journal    = {Invent. Math.},
  year       = {2010},
  volume     = {179},
  number     = {3},
  pages      = {645--681},
  issn       = {0020-9910},
  doi        = {10.1007/s00222-009-0228-0},
  fjournal   = {Inventiones Mathematicae},
  mrclass    = {11F12 (11F66)},
  mrnumber   = {2587342},
  mrreviewer = {Gergely Harcos},
  url        = {https://doi.org/10.1007/s00222-009-0228-0},
}

@Article{BHMM,
  author     = {Blomer, Valentin and Harcos, Gergely and Maga, P\'{e}ter and Mili\'{c}evi\'{c}, Djordje},
  title      = {The sup-norm problem for {$\rm GL(2)$} over number fields},
  journal    = {J. Eur. Math. Soc. (JEMS)},
  year       = {2020},
  volume     = {22},
  number     = {1},
  pages      = {1--53},
  issn       = {1435-9855},
  doi        = {10.4171/jems/916},
  fjournal   = {Journal of the European Mathematical Society (JEMS)},
  mrclass    = {11F72 (11F25 11F55)},
  mrnumber   = {4046009},
  mrreviewer = {Martin Raum},
  url        = {https://doi.org/10.4171/jems/916},
}

@Book{miyake,
  title      = {Modular forms},
  publisher  = {Springer-Verlag, Berlin},
  year       = {1989},
  author     = {Miyake, Toshitsune},
  isbn       = {3-540-50268-8},
  note       = {Translated from the Japanese by Yoshitaka Maeda},
  doi        = {10.1007/3-540-29593-3},
  mrclass    = {11F11 (11F25 11F72)},
  mrnumber   = {1021004},
  mrreviewer = {Harvey Cohn},
  pages      = {x+335},
  url        = {https://doi.org/10.1007/3-540-29593-3},
}

@Book{bump,
  title      = {Automorphic forms and representations},
  publisher  = {Cambridge University Press, Cambridge},
  year       = {1997},
  author     = {Bump, Daniel},
  volume     = {55},
  series     = {Cambridge Studies in Advanced Mathematics},
  isbn       = {0-521-55098-X},
  doi        = {10.1017/CBO9780511609572},
  mrclass    = {11F70 (11F41 11R39 22E50 22E55)},
  mrnumber   = {1431508},
  mrreviewer = {Solomon Friedberg},
  pages      = {xiv+574},
  url        = {https://doi.org/10.1017/CBO9780511609572},
}

@Article{BHMn,
  author     = {Blomer, Valentin and Harcos, Gergely and Maga, P\'{e}ter},
  title      = {Analytic properties of spherical cusp forms on {${\rm GL}(n)$}},
  journal    = {J. Anal. Math.},
  year       = {2020},
  volume     = {140},
  number     = {2},
  pages      = {483--510},
  issn       = {0021-7670},
  doi        = {10.1007/s11854-020-0094-7},
  fjournal   = {Journal d'Analyse Math\'{e}matique},
  mrclass    = {11F72 (11F55 11H06 33E30 43A80)},
  mrnumber   = {4093914},
  mrreviewer = {Rolf Berndt},
  url        = {https://doi.org/10.1007/s11854-020-0094-7},
}

@InCollection{lapid,
  author     = {Lapid, Erez},
  title      = {On the {H}arish-{C}handra {S}chwartz space of {$G(F)\backslash G(\mathbb{A})$}},
  booktitle  = {Automorphic representations and {$L$}-functions},
  publisher  = {Tata Inst. Fund. Res., Mumbai},
  year       = {2013},
  volume     = {22},
  series     = {Tata Inst. Fundam. Res. Stud. Math.},
  pages      = {335--377},
  note       = {With an appendix by Farrell Brumley},
  mrclass    = {11F70 (11F72)},
  mrnumber   = {3156857},
  mrreviewer = {Min Ho Lee},
}

@Article{brumley,
  author     = {Brumley, Farrell},
  title      = {Effective multiplicity one on {${\rm GL}_N$} and narrow zero-free regions for {R}ankin-{S}elberg {$L$}-functions},
  journal    = {Amer. J. Math.},
  year       = {2006},
  volume     = {128},
  number     = {6},
  pages      = {1455--1474},
  issn       = {0002-9327},
  fjournal   = {American Journal of Mathematics},
  mrclass    = {11F70 (11F67 11M20 11M36 22E55)},
  mrnumber   = {2275908},
  mrreviewer = {Daniel Bump},
  url        = {http://muse.jhu.edu/journals/american_journal_of_mathematics/v128/128.6brumley.pdf},
}

@Article{IS,
  author     = {Iwaniec, H. and Sarnak, P.},
  title      = {{$L^\infty$} norms of eigenfunctions of arithmetic surfaces},
  journal    = {Ann. of Math. (2)},
  year       = {1995},
  volume     = {141},
  number     = {2},
  pages      = {301--320},
  issn       = {0003-486X},
  doi        = {10.2307/2118522},
  fjournal   = {Annals of Mathematics. Second Series},
  mrclass    = {11F72 (11F37 58G25 81Q50)},
  mrnumber   = {1324136},
  mrreviewer = {Jens Bolte},
  url        = {https://doi.org/10.2307/2118522},
}

@Article{saha,
  author     = {Saha, Abhishek},
  title      = {Hybrid sup-norm bounds for Maass newforms of powerful level},
  journal    = {Algebra Number Theory},
  year       = {2017},
  volume     = {11},
  number     = {5},
  pages      = {1009--1045},
  issn       = {1937-0652},
  doi        = {10.2140/ant.2017.11.1009},
  fjournal   = {Algebra \& Number Theory},
  mrclass    = {11F03 (11F41 11F60 11F72 11F85 35P20)},
  mrnumber   = {3671430},
  mrreviewer = {Kazuyuki Hatada},
  url        = {https://doi.org/10.2140/ant.2017.11.1009},
}

@Misc{marshall,
  author    = {Marshall, Simon},
  title     = {Upper bounds for Maass forms on semisimple groups},
  year      = {2014},
  copyright = {arXiv.org perpetual, non-exclusive license},
  doi       = {https://doi.org/10.48550/arxiv.1405.7033},
  publisher = {arXiv},
  url       = {https://arxiv.org/abs/1405.7033},
}

@Misc{hu,
  author    = {Hu, Yueke},
  title     = {Sup norm on $\text{PGL}_n$ in depth aspect},
  year      = {2018},
  copyright = {arXiv.org perpetual, non-exclusive license},
  doi       = {https://doi.org/10.48550/arxiv.1809.00617},
  publisher = {arXiv},
  url       = {https://arxiv.org/abs/1809.00617},
}

@Article{toma,
  author   = {Toma, Radu},
  title    = {Hybrid bounds for the sup-norm of automorphic forms in higher rank},
  journal  = {Trans. Amer. Math. Soc.},
  year     = {2023},
  volume   = {376},
  number   = {8},
  pages    = {5573--5600},
  issn     = {0002-9947,1088-6850},
  doi      = {10.1090/tran/8921},
  fjournal = {Transactions of the American Mathematical Society},
  mrclass  = {11F55 (11D45 11F72 11R52)},
  mrnumber = {4630754},
  url      = {https://doi.org/10.1090/tran/8921},
}

@Article{huang,
  author     = {Huang, Bingrong},
  title      = {Sup-norm and nodal domains of dihedral {M}aass forms},
  journal    = {Comm. Math. Phys.},
  year       = {2019},
  volume     = {371},
  number     = {3},
  pages      = {1261--1282},
  issn       = {0010-3616},
  doi        = {10.1007/s00220-019-03335-5},
  fjournal   = {Communications in Mathematical Physics},
  mrclass    = {11F72 (58J50)},
  mrnumber   = {4029832},
  mrreviewer = {Jean Raimbault},
  url        = {https://doi.org/10.1007/s00220-019-03335-5},
}

@Misc{assing,
  author    = {Assing, Edgar},
  title     = {On sup-norm bounds part I: ramified Maaß newforms over number fields},
  year      = {2017},
  copyright = {arXiv.org perpetual, non-exclusive license},
  doi       = {https://doi.org/10.48550/arxiv.1710.00362},
  publisher = {arXiv},
  url       = {https://arxiv.org/abs/1710.00362},
}

@Article{young,
  author     = {Young, Matthew P.},
  title      = {A note on the sup norm of {E}isenstein series},
  journal    = {Q. J. Math.},
  year       = {2018},
  volume     = {69},
  number     = {4},
  pages      = {1151--1161},
  issn       = {0033-5606},
  doi        = {10.1093/qmath/hay019},
  fjournal   = {The Quarterly Journal of Mathematics},
  mrclass    = {11F03 (11F72)},
  mrnumber   = {3908697},
  mrreviewer = {Joshua S. Friedman},
  url        = {https://doi.org/10.1093/qmath/hay019},
}

@Article{AL,
  author     = {Atkin, A. O. L. and Lehner, J.},
  title      = {Hecke operators on {$\Gamma _{0}(m)$}},
  journal    = {Math. Ann.},
  year       = {1970},
  volume     = {185},
  pages      = {134--160},
  issn       = {0025-5831},
  doi        = {10.1007/BF01359701},
  fjournal   = {Mathematische Annalen},
  mrclass    = {10.20},
  mrnumber   = {268123},
  mrreviewer = {R. A. Rankin},
  url        = {https://doi.org/10.1007/BF01359701},
}

@Article{KM,
  author     = {Kenku, M. A. and Momose, Fumiyuki},
  title      = {Automorphism groups of the modular curves {$X_0(N)$}},
  journal    = {Compositio Math.},
  year       = {1988},
  volume     = {65},
  number     = {1},
  pages      = {51--80},
  issn       = {0010-437X},
  fjournal   = {Compositio Mathematica},
  mrclass    = {14G25 (11G30 14K15)},
  mrnumber   = {930147},
  mrreviewer = {Kenneth A. Ribet},
  url        = {http://www.numdam.org/item?id=CM_1988__65_1_51_0},
}

@Article{mcd,
  author     = {McDonald, B. R.},
  title      = {Automorphisms of {${\rm GL}_{n}(R)$}},
  journal    = {Trans. Amer. Math. Soc.},
  year       = {1978},
  volume     = {246},
  pages      = {155--171},
  issn       = {0002-9947},
  doi        = {10.2307/1997969},
  fjournal   = {Transactions of the American Mathematical Society},
  mrclass    = {20G99},
  mrnumber   = {515534},
  mrreviewer = {E. F. Robertson},
  url        = {https://doi.org/10.2307/1997969},
}

@Book{IK,
  title      = {Analytic number theory},
  publisher  = {American Mathematical Society, Providence, RI},
  year       = {2004},
  author     = {Iwaniec, Henryk and Kowalski, Emmanuel},
  volume     = {53},
  series     = {American Mathematical Society Colloquium Publications},
  isbn       = {0-8218-3633-1},
  doi        = {10.1090/coll/053},
  mrclass    = {11-02 (11Fxx 11Lxx 11Mxx 11Nxx)},
  mrnumber   = {2061214},
  mrreviewer = {K. Soundararajan},
  pages      = {xii+615},
  url        = {https://doi.org/10.1090/coll/053},
}

@Article{venk,
  author     = {Venkatesh, Akshay},
  title      = {Large sieve inequalities for {${\rm GL}(n)$}-forms in the conductor aspect},
  journal    = {Adv. Math.},
  year       = {2006},
  volume     = {200},
  number     = {2},
  pages      = {336--356},
  issn       = {0001-8708,1090-2082},
  doi        = {10.1016/j.aim.2005.11.001},
  fjournal   = {Advances in Mathematics},
  mrclass    = {11F55 (11F30 11F66 11N36)},
  mrnumber   = {2200849},
  mrreviewer = {Gergely\ Harcos},
  url        = {https://doi.org/10.1016/j.aim.2005.11.001},
}

@Article{HB,
  author     = {Heath-Brown, D. R.},
  title      = {The density of rational points on curves and surfaces},
  journal    = {Ann. of Math. (2)},
  year       = {2002},
  volume     = {155},
  number     = {2},
  pages      = {553--595},
  issn       = {0003-486X,1939-8980},
  doi        = {10.2307/3062125},
  fjournal   = {Annals of Mathematics. Second Series},
  mrclass    = {11G35 (11G50 14G05 14G40)},
  mrnumber   = {1906595},
  mrreviewer = {Carlo\ Gasbarri},
  url        = {https://doi.org/10.2307/3062125},
}

@Article{schinzel,
  author     = {Schinzel, A.},
  title      = {Reducibility of polynomials of the form {$f(x)-g(y)$}},
  journal    = {Colloq. Math.},
  year       = {1967},
  volume     = {18},
  pages      = {213--218},
  issn       = {0010-1354,1730-6302},
  doi        = {10.4064/cm-18-1-213-218},
  fjournal   = {Colloquium Mathematicum},
  mrclass    = {10.76},
  mrnumber   = {220703},
  mrreviewer = {D.\ J.\ Lewis},
  url        = {https://doi.org/10.4064/cm-18-1-213-218},
}

@Book{newman,
  title      = {Integral matrices},
  publisher  = {Academic Press, New York-London},
  year       = {1972},
  author     = {Newman, Morris},
  volume     = {Vol. 45},
  series     = {Pure and Applied Mathematics},
  mrclass    = {15A33},
  mrnumber   = {340283},
  mrreviewer = {B.\ M.\ Stewart},
  pages      = {xvii+224},
}

@Article{BP,
  author     = {Blomer, Valentin and Pohl, Anke},
  title      = {The sup-norm problem on the {S}iegel modular space of rank two},
  journal    = {Amer. J. Math.},
  year       = {2016},
  volume     = {138},
  number     = {4},
  pages      = {999--1027},
  issn       = {0002-9327,1080-6377},
  doi        = {10.1353/ajm.2016.0032},
  fjournal   = {American Journal of Mathematics},
  mrclass    = {58J50 (11F46 22E46 32N99)},
  mrnumber   = {3538149},
  mrreviewer = {Benjamin\ Linowitz},
  url        = {https://doi.org/10.1353/ajm.2016.0032},
}

@Article{borel-dens,
  author     = {Borel, Armand},
  journal    = {J. Reine Angew. Math.},
  title      = {Density and maximality of arithmetic subgroups},
  year       = {1966},
  issn       = {0075-4102,1435-5345},
  pages      = {78--89},
  volume     = {224},
  doi        = {10.1515/crll.1966.224.78},
  fjournal   = {Journal f\"{u}r die Reine und Angewandte Mathematik. [Crelle's Journal]},
  mrclass    = {14.50 (20.60)},
  mrnumber   = {205999},
  mrreviewer = {H.\ Popp},
  url        = {https://doi.org/10.1515/crll.1966.224.78},
}

@Article{gillman,
  author     = {Gillman, Nate},
  journal    = {J. Number Theory},
  title      = {Explicit subconvexity savings for sup-norms of cusp forms on {${\rm PGL}_n(\mathbb{R})$}},
  year       = {2020},
  issn       = {0022-314X,1096-1658},
  pages      = {46--61},
  volume     = {206},
  doi        = {10.1016/j.jnt.2019.06.002},
  fjournal   = {Journal of Number Theory},
  mrclass    = {11F55},
  mrnumber   = {4013163},
  mrreviewer = {Yuk-Kam\ Lau},
  url        = {https://doi.org/10.1016/j.jnt.2019.06.002},
}

@Article{BT,
  author     = {Brumley, Farrell and Templier, Nicolas},
  title      = {Large values of cusp forms on {$\operatorname{GL}_n$}},
  journal    = {Selecta Math. (N.S.)},
  year       = {2020},
  volume     = {26},
  number     = {4},
  pages      = {Paper No. 63, 71},
  issn       = {1022-1824,1420-9020},
  doi        = {10.1007/s00029-020-00589-z},
  fjournal   = {Selecta Mathematica. New Series},
  mrclass    = {11F70 (58K55)},
  mrnumber   = {4150475},
  mrreviewer = {\Dbar \cftil{o} Ng\d{o}c Di\cfudot{e}p},
  url        = {https://doi.org/10.1007/s00029-020-00589-z},
}

@Misc{dalal-gerbelli,
  author        = {Rahul Dalal and Mathilde Gerbelli-Gauthier},
  title         = {Root Number Equidistribution for Self-Dual Automorphic Representations on $GL_N$},
  year          = {2025},
  archiveprefix = {arXiv},
  eprint        = {2410.01976},
  primaryclass  = {math.NT},
  url           = {https://arxiv.org/abs/2410.01976},
}

@Article{BTZ,
  author     = {Brumley, Farrell and Thorner, Jesse and Zaman, Asif},
  journal    = {J. Eur. Math. Soc. (JEMS)},
  title      = {Zeros of {R}ankin-{S}elberg {$L$}-functions at the edge of the critical strip},
  year       = {2022},
  issn       = {1435-9855,1435-9863},
  note       = {With an appendix by Colin J. Bushnell and Guy Henniart},
  number     = {5},
  pages      = {1471--1541},
  volume     = {24},
  doi        = {10.4171/jems/1134},
  fjournal   = {Journal of the European Mathematical Society (JEMS)},
  mrclass    = {11F66 (11F67)},
  mrnumber   = {4404782},
  mrreviewer = {D.\ R.\ Heath-Brown},
  url        = {https://doi.org/10.4171/jems/1134},
}

\end{document}